\documentclass[12pt]{article}
\usepackage{amssymb,amsmath}
\newtheorem{theorem}{Theorem}[section]
\newtheorem{lemma}{Lemma}[section]

\newtheorem{corollary}{Corollary}[section]
\newtheorem{definition}{Definition}[section]
\newtheorem{example}{Example}[section]

\newcommand{\delp}{\delta \pi}
\newcommand{\w}{\wedge}
\newcommand{\n}{\notag}
\newcommand{\im}{\mbox{i}}
\newcommand{\hook}{\hookrightarrow}
\newcommand{\tb}{\textbf}

\begin{document}

\sloppy

\begin{center}

\LARGE{1-rigidity of CR submanifolds in spheres}

\vspace{1pc}

\large{Sung Ho Wang \\
 Department of Mathematics\\
 Kias \\
 Seoul, Corea 130-722 \\
\texttt{shw@kias.re.kr}}
\end{center}

\vspace{1pc}

\noindent \textbf{Abstract}$\;\,$
We propose a unified computational framework for the problem of deformation
and rigidity of submanifolds in a homogeneous space under geometric constraint.
A notion of 1-rigidity of a submanifold under admissible deformations is introduced.
It means every admissible deformation of the submanifold osculates
a one parameter family of motions up to 1st order.

We implement this idea to the question of rigidity of CR submanifolds in spheres.
A class of submanifolds called Bochner rigid submanifolds are shown to be 1-rigid
under type preserving CR deformations. 1-rigidity is then extended to
a rigid neighborhood theorem, which roughly states that
if a CR submanifold $\, M$ is Bochner rigid, then any pair of mutually CR equivalent
CR submanifolds that are sufficiently close to $\, M$ are congruent
by an automorphism of the sphere.

A local characterization of Whitney submanifold is obtained, which is an example
of a CR submanifold that is not 1-rigid. As a by product, we give a simple
characterization of the proper holomorphic maps from the unit ball
$\, \mathbb{B}^{n+1}$ to $\, \mathbb{B}^{2n+1}$.

\vspace{1pc}


\noindent \textbf{MS classification:} 53B25, 32V30

\thispagestyle{empty}

\newpage
\setcounter{page}{1}

\begin{center}
\LARGE{1-rigidity of CR submanifolds in spheres}

\vspace{1pc}

\large{Sung Ho Wang}
\end{center}

\vspace{1pc}

1. Deformation

2. 1-rigidity

3. CR submanifolds

\quad \quad 3.1 Fundamental forms

\quad \quad 3.2 Bochner rigid submanifolds

4. Rigidity of CR submanifolds

\quad \quad 4.1 CR 1-rigidity

\quad \quad 4.2 Local rigidity

5. Whitney submanifold

6. Proof of Example \ref{brig}

\begin{center}
\textbf{\large{Introduction}}
\end{center}

In the study of submanifolds in a homogeneous space $\, X = \,G/P$
of a Lie group $\, G$, the method of moving frames is both a
unifying concept and an effective tool, which is the version of
the method of equivalence applied to submanifold geometry.
Let $\, \phi$ be the left invariant, Lie algebra $\, \mathfrak{g}$-valued
Maurer-Cartan form of $\, G$. The local equivalence problem for a
submanifold  $\,f: \, M \hook X$ in many interesting cases is solved
on a canonical adapted subbundle $\, E_f: \, B_f \hook G$ together with
$\, \pi = E_f^*\, \phi$ in such a way that the complete set of local
invariants is generated by the coefficients of $\, \pi$ [Ga]. The
pair $( B_f, \, \pi )$ measures in a sense how the submanifold $\,
M$ deviates from a flat model submanifold. The method was
systematically exploited and applied by Cartan himself, and by
Chern in various geometric problems, [Ch] and the references
therein.

A submanifold $\, M$ often inherits a geometric structure from the
ambient space $\, X$, e.g., metric, conformal structure, CR
structure. A general problem of interest is the deformation of a
submanifold preserving the induced geometric structure and the
rigidity phenomena thereof. There exist a wealth of works on the
subject ranging from local differential algebraic analysis of the
prolonged Gau\ss $\,$ equations in the isometric embedding problem
[BBG], to symbol analysis via Lie algebra cohomology in
Griffiths-Harris rigidity of Hermitian symmetric spaces [La][HY]. So
far, however, there appears to be little work which may serve as a
common conceptual ground for these problems.

The purpose of this paper is to propose a unified conceptual perspective
and, more importantly, a computational framework for the problem of
deformation and rigidity of submanifolds in a homogeneous space under geometric
constraint. We introduce a notion of 1-rigidity of a submanifold
under admissible deformations, which can be considered as a
geometric definition of infinitesimal rigidity in classical
differential geometry. The idea is to apply the method of
equivalence to deformation of submanifolds. 1-rigidity means
\emph{every admissible deformation is equivalent to a one parameter family of
motions up to $\,1$st order}. In essence, we differentiate the
method of moving frames once. This approach to rigidity problems may
provide an alternative basis for the existing works
[KT][BBG][La][HY][Ha][ChH].

We implement this idea to the question of rigidity of
CR submanifolds in a sphere $\,  \Sigma^m = \partial \, \mathbb{B}^{m+1}$,
where $\, \mathbb{B}^{m+1}$ is the unit ball in $\mathbb{C}^{m+1}$.
The main result is that
\emph{a class of submanifolds called Bochner rigid
submanifolds are 1-rigid under type preserving CR deformations},
\tb{Theorem \ref{thm3}}. A CR submanifold in $\, \Sigma^m$ is
Bochner rigid if its fundamental forms are rigid in an algebraic
sense, \tb{Definition \ref{bordef}}. 1-rigidity is then extended to
the following local rigidity or rigid neighborhood theorem, \tb{Theorem \ref{lorig}}:
\emph{Let $\,M \subset \Sigma^m$ be a Bochner rigid CR submanifold.
Then any pair of mutually CR equivalent
CR submanifolds that are sufficiently close to $\, M$
are congruent by an automorphism of $\, \Sigma^m$.}
This generalizes the recent work [EHZ] where it was proved, among other things,
that a nondegenerate CR submanifold $\, M^n \subset \Sigma^{n+r}$ of
CR dimension $\, n$ is CR-rigid when $\, r \leq \frac{n}{2}$.

The proof of \tb{Theorem \ref{thm3}} is essentially equivalent to
showing that the linearized Maurer-Cartan equation
\begin{equation}\label{limc}
- \, d (\delp) = \pi \w \, \delp + \, \delp \w \, \pi
\end{equation}
has $\, \delp=0$ as its only solution up to certain group action
(\ref{gpac}). The computation involved is algorithmic as in the
moving frame method, and it proceeds as follows. The condition of
admissible deformation, type preserving CR deformation, imposes a
set of initial conditions on $\, \delp$. Repeated applications of
the linearized Maurer-Cartan equation (\ref{limc}) then give rise to
a sequence of compatibility conditions on $\, \delp$. The Bochner
rigid assumption on the CR submanifold $\, M$ implies that any two
adjacent sequence of compatibility conditions are tightly related,
thus ensuring the propagation of the sequence of compatibility
conditions to force $\, \delp = 0$.

In \tb{Section 1}, the equation of deformation under geometric
constraint is introduced by expanding the Maurer-Cartan equation
of one parameter family of immersions,
for which (\ref{limc}) is the first set of compatibility conditions.
This sequence of equations lead to the notion of 1-rigidity
under admissible deformations in \tb{Section 2}.
Certain group action naturally arises, which plays
an important role analogous to absorption in the method of
equivalence. In \tb{Section 3}, we set up the basic structure
equations for CR submanifolds in spheres, which are then refined for
the class of Bochner rigid submanifolds. The main theorem \tb{Theorem \ref{thm3}}
is proved in \tb{Section 4}, and then extended to local rigidity
theorem \tb{Theorem \ref{lorig}}. The algorithmic computation  for
CR 1-rigidity is carried out in detail in the proof of \tb{Theorem
\ref{thm3}}. In \tb{Section 5}, we present a local characterization
that a CR flat submanifold $\, \Sigma^n \hook \Sigma^{2n}$, $\, n \geq 2$,
is a part of either a linear embedding or a Whitney submanifold, \tb{Theorem \ref{whit}}
The structure equation of a Whitney submanifold together with
Cartan's generalization of Lie's third fundamental theorem implies
that a Whitney submanifold is CR deformable in exactly one direction,
thus providing an example of a CR submanifold which is not 1-rigid.
\tb{Theorem \ref{whit}} also gives a simple characterization of
the proper holomorphic maps from the unit ball
$\, \mathbb{B}^{n+1}$ to $\, \mathbb{B}^{2n+1}$, \tb{Corollary \ref{whitco}}.
\tb{Section 6} is devoted to an algebraic proof that the canonical $\, S^1$-bundles
over Pl\"ucker embeddings of Grassmann manifolds are Bochner rigid.

Our original motivation for the present work is to provide a geometric description of
the complete system for CR maps in [Ha],
the analysis of prolonged Gau\ss $\,$ equations for isometric embeddings in [BBG],
and the Griffiths-Harris rigidity of Hermitian symmetric spaces in [La][HY].
In order to apply our idea to isometric embedding problem,
the higher order terms in the deformation equation (\ref{ddeforme})
should be taken into account.
Application to jet rigidity of certain homogeneous varieties will be reported
elsewhere.

We shall agree that the expression "differentiating A mod B" would
mean "differentiating A and considering mod B".
CR dimension of a CR manifold is by definition the dimension of the contact hyperplane
fields as a complex vector space, \tb{Definition \ref{cr}}. We assume throughout this
paper the CR submanifolds under consideration are of CR dimension $\, \geq 2$.

\section{Deformation}

Let $\, X = G\, / \,P$ be a homogeneous space of a Lie group $\, G$,
and consider a submanifold $\, M \subset X$ defined by an immersion
\begin{equation}\label{sub}
f: \, M \; \hookrightarrow \; X.
\end{equation}
In many geometric situations, the standard reduction process of
moving  frame method gives rise to an adapted subbundle $\,
E_f: B_f \; \hookrightarrow \; G$ with a structure group $\, H
\subset P$ in a canonical way [Ga][Ch].

\begin{picture}(300,80)(-37,0)
\put(93,50){$ f^*G \supset B_f \quad \hookrightarrow \quad G$}
\put(133,30){$\downarrow  H \quad \quad \;  \; \;\downarrow P $}
\put(130,10) {$M \; \quad \hookrightarrow \quad X$}
\put(165,20){$f$}
\put(160,60){$E_f $}
\end{picture}

Such $\, B_f$ captures the geometry of $\, f$ in the following
sense.  Let $\, \mathfrak{g}$ be the Lie algebra of $\, G$, and
let $\, \phi$ be the $\, \mathfrak{g}$-valued left invariant Maurer-Cartan form of $\, G$:
Let $\, f_1$, $\, f_2$ be two immersions with the induced $\,
H$-bundles $\, B_{f_1}$, $\, B_{f_2}$ respectively. $\, f_1$ and $\,
f_2$ are congruent up to a left motion by an element of $\, G$ iff
there exists a $\, H$-bundle isomorphism of the pair $\, ( B_{f_1},
\, E_{f_1}^* \phi)$ and $\, ( B_{f_2}, \, E_{f_2}^* \phi)$
[Br1][Gr].

We wish to describe in this section the geometry of the deformation
of a submanifold in this setting. Assume for definiteness $\, g: G
\hookrightarrow \mbox{GL}(N, \, \mathbb{R})$ or $\, \mbox{GL}(N, \,
\mathbb{C})$ for an integer $\, N$. Then in particular $\, \phi=
g^{-1} \, dg$, and $\, \phi$ satisfies  Maurer-Cartan structure
equation
\[ -\,d \phi =  \phi \w \phi. \]

Consider a deformation $\, f_t$ of $\, f_0=f$ for  $\, t\in(
-\epsilon, \, \epsilon)$, and let $\, E_t: B_t \; \hookrightarrow \;
G$ be the associated adapted $\, H$-bundle.

\begin{picture}(300,80)(-37,0)
\put(93,50){$ f_t^*G \supset B_t \quad \hookrightarrow \quad G$}
\put(134,30){$\downarrow  H \quad \; \;  \quad \downarrow  P$}
\put(130,10) {$M  \quad \hookrightarrow \quad X$}
\put(160,20){$f_t$}
\put(160,60){$E_t$}
\end{picture}

\noindent A general element $\, E_t = ( e_1(t), \, e_2(t), \, ... \,
, \, e_N(t))$ may be written as
\begin{equation}
e_A(t) = e_B(0) \, g^B_A(t), \, \; A, \, B \, = 1, \, ... \, , N. \n
\end{equation}
Equivalently in matrix form
\begin{equation}
E_t = E_0 \, g_t, \n
\end{equation}
where $\, g_t = (\, g^B_A(t) \, )$ is the unique $\, G$-valued
function on $\, B_f \times \, (-\epsilon, \, \epsilon)$ that
describes the deformation of the bundle  $\, B_f$. Let us expand $\,
g_t$ (formally) with respect to $\, t$
\begin{equation}
g_t = \exp(\sum_{k=1}^{\infty} \frac{t^k}{k!}\, U_k) \n
\end{equation}
for a sequence of $\, \mathfrak{g}$-valued functions $\{\, U_k \,
\}$  on $\, B_f$. Then
\begin{align}
dE_t &= dE_0 \, g_t + E_0 \, dg_t \notag \\
&= E_0 \, \pi_0 \, g_t + E_0 \, dg_t \notag \\
&= E_t \, ( g_t^{-1} \, \pi_0 \, g_t + g_t^{-1} \, dg_t) \notag \\
&= E_t \, \pi_t   \n
\end{align}
where $\, \pi_t = E_t^{-1} \, dE_t = E_t^* \, \phi$. The expression
for $\, \pi_t = \, \sum_{k=1}^{\infty} \, \frac{t^k}{k!}\, \pi_k$
can be computed as follows.
\begin{align}
\pi_1 &= dU_1 + \pi_0 \, U_1 - U_1 \, \pi_0 \label{deforme} \\
\pi_2 &= dU_2 + \pi_0 \, U_2 - U_2 \, \pi_0  + \pi_1 \, U_1 - U_1 \, \pi_1 \notag \\
&\vdots \notag \\
\pi_k &\equiv dU_k + \pi_0 \, U_k - U_k \, \pi_0 \; \; \mod \; \{ \,
\pi_{k-1}, \, \pi_{k-2}, \, ... \, \pi_1 \}\notag \\
&\vdots \notag
\end{align}

On the other hand, the deformations we are interested in are not
arbitrary ones but those deformations that satisfy certain geometric
constraint. Since the induced Maurer-Cartan form $\, \pi_t$ captures
the geometry of the immersion $\, f_t$, it is reasonable to define
the condition of admissible deformation in terms of \emph{a finite
set of ordinary differential relations among the coefficients of $\,
\pi_t$}. These relations are in turn expressed as a sequence
of (algebraic) relations among the coefficients of
$\, \{ \, \pi_1, \, \pi_2, \, ... \, \}$.

\emph{Remark.}$\;$ There are situations where the condition of
admissible deformation is given in a form which is more general than
a sequence of algebraic relations among the coefficients of
$\, \{ \, \pi_1, \, \pi_2, \, ... \, \}$, e.g., the jet rigidity of
homogeneous varieties in projective spaces [LM].
The condition of admissible deformation in this case
may involve a sequence of differential relations among the coefficients of
$\, \{ \, \pi_1, \, \pi_2, \, ... \, \}$.

\begin{definition}
Let $\, f: M \hookrightarrow X$ be a submanifold in a homogeneous
space $\, X = G/P$. Suppose a condition of admissible deformation is
given as a finite set of ordinary differential relations among the
coefficients of $\, \pi_t$. The \emph{equation of deformation} at
$\, f_0=f$ is the sequence of first order equations \eqref{deforme}
for  a sequence of $\, \mathfrak{g}$-valued functions $\{ \, U_1, \,
U_2, \, ... \, \}$ on the adapted $\, H$-bundle $\, B_f$ with $\{ \,
\pi_1, \, \pi_2, \, ... \, \}$ satisfying the sequence of relations
that correspond to the given set of differential relations on $\,
\pi_t$.
\end{definition}
The idea of this set up is to lift the problem of deformation to
Lie group $\, G$ where a uniform treatment is available through
Maurer-Cartan form and its structure equation.

Note by expanding the structure equation
\begin{equation}
- \, d\pi_t = \, \pi_t \, \wedge \pi_t \n
\end{equation}
with respect to $\, t$, we obtain a sequence of differential
equations for $\{ \, \pi_1, \, \pi_2, \, ... \, \}$.
\begin{align}
-\, d \pi_1 &=  \pi_0 \wedge \, \pi_1 + \pi_1 \wedge \,  \pi_0, \label{ddeforme}\\
-\, d \pi_2 &=  \pi_0 \wedge \, \pi_2 + \pi_2 \wedge \, \pi_0  + 2 \, \pi_1 \wedge \,\pi_1, \notag \\
&\vdots \notag \\
-\, d \pi_k &= \sum_{i+j=k}  \, \frac{k!}{i!j!}\, \pi_i \wedge \, \pi_j, \notag \\
&\vdots \notag
\end{align}
These are in fact the set of compatibility conditions obtained by
differentiating (\ref{deforme}) once.

We now introduce a certain \emph{group action} on the space of
deformations that is tangent to Cauchy characteristics of the
deformation equation, which plays an important role in our treatment
of rigidity questions.

Assume a deformation $\, E_t$ is given. Suppose $\, h_t$ is a  $\,
H$-valued function on $\, B_f \times \, (-\epsilon, \, \epsilon)$,
and consider  a new deformation
\begin{equation}
\hat{E_t} = E_t \, h_t. \n
\end{equation}
Let $\, \it{\Pi}: \, G \to \, X = G/P$ be the projection map. Since
$\, H \subset P$, $\, \it{\Pi}(\hat{E_t}) =  \it{\Pi}(E_t) = f_t$ describes
the same deformation of the immersion $\, f$. This action of
$\, C^{\infty}(B_f \times \, (-\epsilon, \, \epsilon), \, H)$ simply
rotates along the fiber of $\, B_t \to M$.

The nature of this ambiguity is similar to the one in the reduction
process of moving frame method. In order to see this group action
explicitly, let us write $\, h_t = \exp(\sum_{k=1}^{\infty}
\frac{t^k}{k!}\, V_k)$ for a sequence of $\, \mathfrak{h}$-valued
functions $\, \{ \, V_1, \, V_2, \, ... \, \}$ on $\, B_f$, where
$\, \mathfrak{h}$ is the Lie algebra of $\, H$. Then $\, V_k$'s
\emph{acts} on $\, \pi_k$'s by translation as follows.
\begin{align}
\Delta \pi_1 &=      \pi_0 \, V_1 - V_1 \, \pi_0 \label{gpac}\\
\Delta \pi_2 &\equiv \pi_0 \, V_2 - V_2 \, \pi_0 \; \; \mod \; \{ \, V_1 \}   \notag \\
&\vdots \notag \\
\Delta \pi_k &\equiv \pi_0 \, V_k - V_k \, \pi_0 \; \; \mod \; \{ \, V_{k-1}, \, V_{k-2},
\, ... \, V_1 \}\notag \\
&\vdots  \notag
\end{align}
Here $\, \Delta \pi_k$ stands for the contribution from  $\, \{ \,
V_1, \, V_2, \, ... \, \}$. We will  utilize this extra degree of
freedom to normalize certain coefficients in the rigidity related
computations. This is analogous to the process of \emph{absorption}
in the method of equivalence [Ga].

\section{1-rigidity}

The notion of 1-rigidity under admissible deformations of a
submanifold in a homogenous space is introduced in this section.
1-rigidity can be considered as a geometric and unifying definition
of infinitesimal rigidity in classical differential geometry. We
continue to use the notations adopted in \textbf{Section 1}.

Suppose in the deformation equation (\ref{deforme})
\begin{equation}
\pi_i=0  \quad \, \mbox{for} \; 1 \, \leq i \, \leq k, \n
\end{equation}
or equivalently
\begin{equation}
dU_i + \pi_0 \, U_i - U_i \, \pi_0 =0 \, \quad \mbox{for}
\; 1 \, \leq i \, \leq k. \n
\end{equation}
Then it is easy to check $\, U_i = E_0^{-1} \, a_i \, E_0$ for a
constant $\, a_i \in \mathfrak{g}$, $1 \; \leq i \; \leq k$, and
\begin{align}
E_t &= E_0 \, g_t = E_0 \, \exp(\frac{t^i}{i!}\, U_i) \notag \\
&\equiv  \, \exp(\sum_{i=1}^k \frac{t^i}{i!}\, a_i) \, E_0
\quad \mod \; t^{k+1}.\n
\end{align}
The vanishing of the first $k$ $\, \pi_i$'s thus means that  the
deformation osculates a one parameter family of motions by $\, G$ up
to order $k$.

\begin{definition}
Suppose a condition of admissible deformations is given as a finite
set of ordinary differential relations among the coefficients of $\,
\pi_t$. A submanifold in a homogeneous space is \emph{1-rigid} under
admissible deformations if every solution to the deformation equation \eqref{deforme}
necessarily has $\, \pi_1=0$ modulo the group action \eqref{gpac}.
\end{definition}

The original deformation equation (\ref{deforme}) is however rather difficult
to use in practice to verify 1-rigidity of a submanifold.
Instead, note $\, \pi_1$ satisfies the compatibility equation (\ref{ddeforme}),
\[ -d\pi_1 = \pi_0 \wedge \pi_1 + \pi_1 \wedge \pi_0.\]
For all practical purposes, it is this equation modulo the group action (\ref{gpac})
that we will use to test for 1-rigidity.
The process of simplification of this sort, that is removing
nongeometric terms by differentiating once, has been capitalized in [BG]
in their study of characteristic cohomology.

As it is the case with the method of moving frames, this definition is
best explained and understood in practice. We consider in this section
as an example a submanifold in the homogeneous space of
conformal sphere $\, S^m$, and its conformal deformation. In the
course of computation, a set of criteria for conformal 1-rigidity of
a submanifold naturally emerges. This example also serves as a guide
through the computation for CR 1-rigidity in \textbf{Section 4}. For
general reference for conformal geometry, [Br2].

Let $\, \mathbb{R}^{m+1,1}$ be given a metric of signature
$(m+1,1)$. Let $\, G \, =$ SO$_0(m+1,1)$ be the identity component
of the group of  linear transformations that preserve the metric.
Let $\, X \, = \, S^m$ be the set of \emph{positive} null rays
through the origin. It is well known that SO$_0(m+1,1)$ acts
transitively on $\, S^m$, and that $\, S^m$ inherits an
SO$_0(m+1,1)$ invariant conformal structure.

Let $\, f: \, M \, \hookrightarrow X$ be an $n$-dimensional
submanifold. Upon an appropriate choice of basis for $\, \mathbb{R}^{m+1,1}$,
we may arrange so that the induced Maurer-Cartan form $\, \pi$ on
a \emph{1-adapted} $\, H$-bundle $\,
B_f$ takes the following form, [Br2] for details on 1-adapted bundle $\, B_f$.
\begin{equation}\label{const}
\pi=
\begin{pmatrix}
\pi_0^0 & \pi_j^0 & \pi_b^0 & 0           \\
\pi_0^i & \pi_j^i & \pi_b^i & \pi_{m+1}^i \\
0       & \pi_j^a & \pi_b^a & \pi_{m+1}^a \\
0       & \pi_j^{m+1} & 0   & \pi_{m+1}^{m+1}
\end{pmatrix}
\end{equation}
where $\, \pi_j^i = -\pi_i^j$, $\, \pi_b^a = -\pi_a^b$,
$\,\pi_j^{m+1} = - \pi_0^j$, $\, \pi_j^0 = - \pi_{m+1}^j$, $\,
\pi_b^0 = - \pi_{m+1}^b$, $\, \pi_{m+1}^{m+1}= - \pi_0^0$, and $\, 1
\leq i, \, j \leq n$,  $\, n+1 \leq a, \, b \leq m$. $\, \pi_0^i$'s
are semibasic 1-forms, and the quadratic form $\, \, \pi_0^i \circ \pi_0^i$
represents the induced conformal structure on $\, M$.

$\, \pi$ satisfies the structure equations $\, -\, d\pi= \, \pi \w
\pi$.  Differentiating $\, \pi^a_0 = 0$, we get $\, \pi^a_i \w \pi^i_0 = 0$,
and by Cartan's Lemma
\begin{equation}
\pi_i^a = h^a_{ij} \, \pi_0^j \n
\end{equation}
for a set of coefficients $\, h^a_{ij}=h^a_{ji}$. Using the group
action by $\, \pi^a_{m+1}$, we normalize so that  tr $h^a_{ij} = 0$.
The trace free parts of the second fundamental forms $\, H^a \,  =
\, h^a_{ij} \, \pi_0^i \circ  \pi_0^j$ are conformal invariant of
the submanifold $\, M$. Note that an element of the Lie algebra $\, \mathfrak{h}$
of $\, H$ is of the following form.
\begin{equation}
\begin{pmatrix}
V_0^0 & V_j^0 &  0    & 0         \\
0     & V_j^i &  0    & V_{m+1}^i \\
0       & 0 & V_b^a & 0         \\
0       & 0 & 0   & V_{m+1}^{m+1}
\end{pmatrix}\n
\end{equation}

Suppose we are interested in  1-rigidity under the deformations
that preserve the induced conformal structure on $\, M$. From
(\ref{const}), this is written as an ordinary differential equation
\begin{equation}
\frac{d}{dt} \, \pi_0^i(t) \circ \pi_0^i(t) |_{t=0} = \, \lambda  \, \pi_0^i(0) \circ \pi_0^i(0) \n
\end{equation}
for a scaling factor $\, \lambda$. Using the group action by
$\, V_0^0, \, V_j^i$ components, however, we may in fact
normalize so that
\begin{equation}
\frac{d}{dt} \, \pi_0^i(t)|_{t=0} = 0.\n
\end{equation}
Hence $\, \pi_1 = \delta \pi$ is of the following form.
\begin{equation}\label{confd}
\delta \pi=
\begin{pmatrix}
\delta \pi_0^0 & \delta \pi_j^0 & \delta \pi_b^0 & 0           \\
0 & \delta \pi_j^i & \delta \pi_b^i & \delta \pi_{m+1}^i \\
0       & \delta \pi_j^a & \delta \pi_b^a & \delta \pi_{m+1}^a \\
0       & 0 & 0   & \delta \pi_{m+1}^{m+1}
\end{pmatrix}
\end{equation}
The remaining \emph{group variables} at this stage are
$\, V_b^a$ and $\, V^0_i = - V_{m+1}^i$.

Differentiating $\, \delp^i_0=0, \, \delp^a_0=0$ in (\ref{confd})
using (\ref{deforme}), and by Cartan's lemma, it is easy to check
there exists a set of coefficients $\, p_i, \, p^a_{ij}=p^a_{ji}$
such that
\begin{align}
\delta \pi_0^0 &= p_i \, \pi_0^i, \notag \\
\delta \pi_j^i &= p_i \, \pi_0^j - p_j \, \pi_0^i,\notag \\
\delta \pi_i^a &= p^a_{ij} \, \pi_0^j.\n
\end{align}
Using the group action by  $\, V_j^0$ component as above,  we may
translate $\, p_j =0$. Now differentiating $\, \delta \pi_0^0 =0, \,
\delta \pi_j^i=0$ with these relations, we obtain the following
compatibility conditions for deformation.
\begin{align}
\delta \pi_i^0 &= p_{ij} \, \pi_0^j,   \quad \; p_{ij}=p_{ji}, \notag \\
\pi^a_i \wedge \, &\delta \pi^a_j + \delta \pi^a_i \wedge \pi^a_j -
\pi^i_0 \wedge \, \delta \pi^0_j - \delta \pi^0_i \wedge \pi^j_0=0. \label{confint}
\end{align}

Fix a point $\, p \in M$, and let $\, V^n = T_p M$,  $\, W^{m-n} =
N_p M$ represent the tangent space and the normal space
respectively. Take a conformal basis $\, \{ x^i \}$ for $\, V^*$ and
$\, \{ w_a \}$ for $\, W$. Set
\begin{align}
H = h^a_{ij} \, x^i \, x^j \otimes w_a \in S^2 V^* \otimes W, \notag \\
P = p^a_{ij} \, x^i \, x^j \otimes w_a \in S^2 V^* \otimes W.\n
\end{align}
Let $\, K(V) \subset \bigwedge^2 V^* \otimes \bigwedge^2 V^*$ be
the space of curvature like tensors, which decomposes into
$\, K(V) = \mathcal{W} \oplus Ric$, Weyl curvature tensor and Ricci
tensor. Then the equation (\ref{confint}) is equivalent to
\begin{align}
&\gamma(H, P)^{\mathcal{W}} =0,  \label{weyl1} \\
&\gamma(H, P)^{Ric}  = p_{ij} \, x^i \, x^j \label{weyl2}
\end{align}
where $\, \gamma(H, P)$ is the linearized Gauss map [BBG]. An
element $\, H  \in S^2 V^* \otimes W$ is called \emph{Weyl rigid}
when the space of solutions to (\ref{weyl1}) is $\, \{ \, P = v^a_b
\, h^b_{ij} \, x^i \, x^j \otimes w_a \, | \; \, v^a_b = -\, v^b_a \, \}$.
Note for any $\, Q \in S^2 V^* \otimes W$,
$\, \gamma(Q, P)^{\mathcal{W}} = \gamma(Q_0, P)^{\mathcal{W}}$,
where $\, Q_0$ is the trace free part of $\, Q$.
We mention in passing that when dim$\, W=1$, a traceless
quadratic form $\, H$ is Weyl rigid if it has 3 nonzero eigenvalues.

Assume the trace free part of the second fundamental form $\, H$ of
$\, M$  is Weyl rigid. Then using the group action by $\, V^a_b$
component, we may translate $\, \delta \pi_i^a =0$. It follows $\,
\delta \pi_i^0 =0$ by (\ref{weyl2}), and the only remaining nonzero
elements of $\, \delta \pi$ are $\, \delta \pi^0_b$ and $\, \delta
\pi^a_b$. Differentiating $\, \delta \pi_i^a =0, \, \delta \pi_i^0
=0$, we finally get
\begin{align}
\delta \pi_b^0 \, \wedge \pi_i^b = 0,  \label{confla1} \\
\delta \pi_b^a \, \wedge \pi_i^b= 0.   \label{confla2}
\end{align}

An element $\, H \in  S^2 V^* \otimes W$ is nondegenerate if
the associated map $\, H:  \, S^2 V \to  W$ is surjective.
We first show that for a nondegenerate $\, H$, (\ref{confla2}) implies
$\, \delp_b^a=0$. Let $\, \delp^a_b = q^a_{bi} \, \pi^i$, $q^a_{bi} = - q^b_{ai}$.
Then (\ref{confla2}) is equivalent to
$
\,  q^a_{bi} \, h^b_{jk} = q^a_{bk} \, h^b_{ji} = q^a_{ijk}
$
with $\, q^a_{ijk}$ fully symmetric in lower indices. Multiplying
both sides by $\, h^a_{sp}$ and summing over $\, a$,
\begin{align}
q^a_{ijk} \, h^a_{sp} &= q^a_{bi} \, h^b_{jk} \,  h^a_{sp} \n\\
                      &= -\, q^b_{ai} \, h^b_{jk} \,  h^a_{sp} \n\\
                      &= -\, q^b_{isp} \, h^b_{jk}. \n
\end{align}
But $\, q^a_{ijk}$ is symmetric in lower indices and $\, h^a_{ij}$ is nondegenerate,
hence $\, q^a_{ijk} =0$ and consequently $\, q^a_{bi} =0$.

We next show if $\, H$ is both nondegenerate and Weyl rigid, then
(\ref{confla1}) implies $\, \delp_b^0=0$.
Let $\, \delp^0_b = q_{bi} \, \pi^i$.
Then (\ref{confla1}) is equivalent to
$
\,  q_{bi} \, h^b_{jk} = q_{bk} \, h^b_{ji} = q_{ijk}
$
with $\, q_{ijk}$ fully symmetric in lower indices. For an arbitrary set
of coefficients $\, \{ \, f_i \, \}$, set
\[ Q = Q^a \otimes e_a = (q_{ai} \, f_j + q_{aj} \, f_i) \, x^i  x^j  \otimes e_a
\in S^2V^* \otimes W. \]
Then one easily sees $\, \gamma(H, \, Q)^{\mathcal{W}}=0$.
Since $\, H$ is Weyl rigid, $ q_{ai} \, f_j + q_{aj} \, f_i = v^a_b\, h^b_{ij} $
for  skew symmetric $\, v^a_b = - \, v^b_a$.
Multiplying both sides by $\, h^a_{sp}\, x^i x^j x^s x^p$
and summing over $\, a,  i,  j,  s,  p$, we get
\[ (f_i \, x^i) \, ( q_{jsp} \,  x^j\, x^s\, x^p) =0.\]
Since $\, \{ \, f_i \, \}$ is arbitrary, $\, q_{ijk} =0$, and consequently $\, q_{bi} =0$
for $\, h^a_{ij}$ is nondegenerate.

Thus \emph{a submanifold $\, M \hookrightarrow S^m$ is conformally 1-rigid if it
has a nondegenerate Weyl rigid second fundamental form}.
One can apply the results of [KT] and show that the canonical isometric
embeddings of irreducible compact Hermitian symmetric spaces into spheres
satisfy this criteria, and hence they are conformally 1-rigid.

\section{CR submanifolds}
\subsection{Fundamental forms}\label{Fundaforms}
In this section, we set up the basic structure equations for CR
submanifolds  in spheres. A sequence of invariants called
\emph{fundamental forms} are derived from the structure equations.
For general reference in CR geometry, [ChM][Ja].

Let $\, \mathbb{C}^{m+1,1}$ be the complex vector space with coordinates
$\, z = (\, z^0, \, z^A, \, z^{m+1} \,)$, $\, 1 \, \leq A \, \leq  m$,
and a Hermitian scalar product
\begin{equation}
\langle \, z, \, \bar{z} \, \rangle = z^A \, \bar{z}^A +
\mbox{i} \, (z^0 \, \bar{z}^m - z^m \, \bar{z}^0 ).\n
\end{equation}
Let $\, \Sigma^m$ be the set of equivalence classes up to scale of
null vectors with respect to this product. Let SU$(m+1,1)$ be the
group of unimodular linear transformations that leave the form $\,
\langle \, z, \, \bar{z} \rangle$ invariant. Then SU$(m+1,1)$ acts
transitively on $\, \Sigma^m$, and
\begin{equation}
p: \, \mbox{SU}(m+1,1) \to \Sigma^m = \, \mbox{SU}(m+1,1)/P \n
\end{equation}
for an appropriate  subgroup $\, P$ [ChM].

Explicitly, consider an element $\, Z=(\, Z_0, \, Z_A, \, Z_{m+1}\,)
\in \mbox{SU}(m+1,1)$ as an ordered set of $\, (m+2)$-column vectors
in $\, \mathbb{C}^{m+1,1}$  such that det$(Z)=1$, and that
\begin{equation}\label{prod}
\langle \, Z_A, \, \bar{Z}_B  \, \rangle = \, \delta_{AB}, \quad \langle \,Z_0,
\, \bar{Z}_{m+1} \, \rangle = - \langle \,Z_{m+1}, \, \bar{Z}_0 \, \rangle=\, \mbox{i},
\end{equation}
while all other scalar products are zero. We define $\, p(Z) \, = \,
[Z_0]$, where $\, [Z_0]$ is the equivalence class of null vectors
represented by   $\, Z_0$. The left invariant Maurer-Cartan form $\,
\phi$ of SU$(m+1,1)$ is defined by the equation
\begin{equation}
d\,Z = Z \, \phi,\n
\end{equation}
which is in coordinates
\begin{equation}\label{struct1}
d (\, Z_0, \, Z_A, \, Z_{m+1} \, ) = ( \, Z_0, \, Z_B, \, Z_{m+1} \, ) \,
\begin{pmatrix}
\phi_0^0 & \phi_A^0 & \phi_{m+1}^0 \\
\phi_0^B & \phi_A^B & \phi_{m+1}^B \\
\phi_0^{m+1} & \phi_A^{m+1} & \phi_{m+1}^{m+1}
\end{pmatrix}.
\end{equation}
Coefficients of $\, \phi$ are subject to the relations obtained from
differentiating (\ref{prod}) which are
\begin{align}
\phi_0^0 + \bar{\phi}_{m+1}^{m+1} &=0 \n \\
\phi_0^{m+1}&= \bar{\phi}_0^{m+1}, \quad \phi_{m+1}^{0}= \bar{\phi}_{m+1}^{0} \notag \\
\phi_A^{m+1}&= - \, \mbox{i} \, \bar{\phi}_0^A, \quad \phi_{m+1}^A = \mbox{i} \, \bar{\phi}_A^0 \notag \\
\phi_B^A + \bar{\phi}_A^B &=0 \notag \\
\mbox{tr} \,  \phi &=0,  \n
\end{align}
and $\, \phi$ satisfies the structure equation
\begin{equation}\label{struct2}
- \, d\phi = \phi \wedge \phi.
\end{equation}

It is well known that the SU$(m+1,1)$-invariant CR structure on $\,
\Sigma^m \subset \mathbb{C}P^{m+1}$ as a real hypersurface is
biholomorphically equivalent to the standard  CR structure on
$\, S^{2m+1} = \partial \mathbb{B}^{m+1}$, where
$\, \mathbb{B}^{m+1} \subset \mathbb{C}^{m+1}$ is the unit ball.
The structure equation (\ref{struct1}) shows that for any local section
$\, s: \, \Sigma^m \, \to$ SU$(m+1,1)$, this CR structure is defined by
the hyperplane fields $\, (s^*\phi_0^{m+1})^{\perp} = \mathcal{H}$
and the set of (1,0)-forms $\, \{ \, s^*\phi_0^A \, \}$.

\begin{definition}\label{cr}
Let $\, M$ be a manifold of dimension $\, 2n+1$. A submanifold
defined  by an immersion $\, f: \, M \, \hookrightarrow \, \Sigma^m$
is a \emph{CR submanifold} if $\, f_*T_pM \, \cap \,
\mathcal{H}_{f(p)}$ is a complex subspace of $\, \mathcal{H}_{f(p)}$
of dimension $\, n$ for each $\, p \in M$.
\end{definition}
Note the induced hyperplane fields $\, f_*^{-1}(f_*(TM) \, \cap \,
\mathcal{H})$ is necessarily a contact structure on $\, M$, and thus
$\, M$ has  an induced nondegenerate CR structure.

Let $\, f:\, M \, \hookrightarrow \, \Sigma^m$ be a CR submanifold.
For each $\, p \in M$, we wish to define an associated decomposition
of $\, \mathbb{C}^{m+1,1}$, which is a part of the reduction process
of moving frame method applied to CR submanifold.
Let $\, \Pi: \, \mathbb{C}^{m+1,1}-\{0\} \to \mathbb{C}P^{m+1}$ denote the projection map.

Set $\, V_0 = \, \langle \Pi^{-1}(f(p)) \, \rangle$ and
$\, V_0+V_1 = \, \langle  \Pi^{-1}(f(p)), \, \Pi^{-1}(\partial \, f) \,
\rangle$, where $\, \partial  f$ stands for the holomorphic or (1,0)
derivatives of  $\, f$. Successively define the sequence of
subspaces $\, W_2, \, W_3, \, ...\, , \, W_{\tau}$ with dim $W_l=
r_l$ so that
\begin{equation}
(V_0+V_1) \oplus W_2 \oplus \, ... \,  W_l
= \, \langle \, \Pi^{-1}(f(p)), \, \Pi^{-1}(\partial  f), \, \Pi^{-1}(\partial^2  f), \, ... \,
\, \Pi^{-1}(\partial^l f) \, \rangle
\end{equation}
and $\, n+r_2+ \, ... \, r_{\tau} = m$, where $\, \partial^l   f$
stands for the $l$-th order (1,0) derivatives  of $\, f$.
The set of numbers $(\, r_2, \, r_3, \, ... \, r_{\tau} \, )$ are called
the \emph{type numbers} of $\, f$ at $\, p$, and $\, \tau$ is called the \emph{height}.

\emph{Remark.}$\;$ The structure equation \eqref{struct1},  \eqref{hs},
and  $\, \pi^0_A = \mbox{i}\, \bar{\pi}^A_{m+1}$ show
this orthogonal decomposition, rather than a filtration, is well defined.

\begin{definition}
A CR submanifold is of \emph{constant type} if the type numbers are
constant.
\end{definition}
We assume $\, M \hookrightarrow \Sigma^m$ is a CR submanifold of
constant type from now on.

In terms of the structure equation (\ref{struct1}), we may arrange so that
\begin{align}
V_0 &= \, \langle  \, Z_0 \, \rangle \label{decom} \\
(V_0+V_1) &= \, \langle  \, Z_0  \, \rangle + \langle \,  Z_i \, \rangle_{i=1}^{n} \notag \\
(V_0+V_1) \oplus W_2 &= \, \big{(} \langle  \, Z_0  \, \rangle
+ \langle \,  Z_i \, \rangle_{i=1}^{n} \big{)}
\oplus \, \langle   \, Z_{i_2} \, \rangle_{i_2=n+1}^{n_2} \notag \\
&\vdots \notag \\
(V_0+V_1) \oplus W_2 \oplus \, ... \,  W_l &=  \, \big{(} \langle  \, Z_0  \, \rangle
+ \langle \,  Z_i \, \rangle_{i=1}^{n} \big{)}
 \oplus \, \langle  \, Z_{i_2} \, \rangle_{i_2=n+1}^{n_2} \oplus \, ... \,
\, \langle  \, Z_{i_l} \, \rangle_{i_l=n_{l-1}+1}^{n_{l}} \notag \\
&\vdots \notag \\
(V_0+V_1) \oplus W_2 \oplus \, ... \,  W_{\tau} &=  \, \big{(} \langle  \, Z_0  \, \rangle
+ \langle \,  Z_i \, \rangle_{i=1}^{n} \big{)}
 \oplus \, \langle  \, Z_{i_2} \, \rangle_{i_2=n+1}^{n_2} \oplus \, ... \,
\, \langle  \, Z_{i_{\tau}} \, \rangle_{i_{\tau}=n_{\tau-1}+1}^{n_{\tau}}, \notag
\end{align}
where $\, n_l = n_{l-1}+r_l$ and $\, n_1 = n$.
This decomposition then induces an associated adapted  subbundle $\, B$
\begin{equation}
E: \, f^*(\mbox{SU$(m+1,1)$}) \supset B \, \hook \mbox{SU$(m+1,1)$}  \n
\end{equation}
on which the following weak structure equations hold. Let  $\, \pi = E^*\phi$.
\begin{align}\label{wstruct}
\pi_0^A &= 0 \quad \; \; \; \, \mbox{for} \; \;  n+1 \leq \, A \, \leq \, m,    \\
\pi_i^A &\equiv 0 \quad  \mod \; \; \theta, \; \pi^{0,1}
\quad \mbox{for} \; \;  n_2+1 \leq \, A \, \leq \, m,   \n  \\
\pi_{i_l}^A &\equiv 0 \quad  \mod \; \; \theta, \; \pi^{0,1}
\quad \mbox{for} \; \;  n_{l+1}+1 \leq \, A \, \leq \,   m,
\quad 2 \leq l \leq \tau-2,  \n
\end{align}
where $\, n_{l-1}+1 \leq i_l \leq n_l$ for $\, 2 \leq l \leq \tau$,
and$\mod \pi^{0,1}$
means$\mod \bar{\pi}_0^1, \, \bar{\pi}_0^2, \, ... \, \bar{\pi}_0^n$.
\eqref{wstruct} implies the following equation in matrix form.
\begin{equation}
\pi \equiv
\begin{pmatrix}
*         & *              & *     & *       &       & *   \\
\pi^i_0     & *              & *     & *       & ...   & *   \\
0         & \pi^{i_2}_{i_1}             & *     & *       &       & *   \\
0         & 0              & \pi^{i_3}_{i_2}     & *       &       & *   \\
&               &  \vdots    &      & \ddots      & *   \\
0         & 0              & 0     & 0     &       & *
\end{pmatrix} \quad \mod \pi^{m+1}_0, \pi^{0,1}. \n
\end{equation}

Denote $\, \pi_0^{m+1} = \, \theta$,  $\pi^i=\pi^i_0$.
Differentiating $\, \pi^{A}_0=0$  $\mbox{for} \; \;  n+1 \leq \, A \, \leq \, m$,
we get      $\, \pi^{A}_i \wedge \pi^i +  \pi^{A}_{m+1}
\wedge \theta =0$. By Cartan's lemma, there exists a set of
coefficients $\,h^{A}_{ij}=h^A_{ji}, \, h^{A}_i, \, h^{A}$ such that
\begin{equation}\label{hs}
\begin{pmatrix}
\pi^{A}_i \\
\pi^{A}_{m+1}
\end{pmatrix}
=
\begin{pmatrix}
h^{A}_{ij} & h^{A}_i \\
h^{A}_j  & h^{A}
\end{pmatrix}
\begin{pmatrix}
\pi^j \\
\theta
\end{pmatrix}.
\end{equation}
Note by definition of $\, W_2$, $\, h^A_{ij} = 0\;$ for $\, A \geq n_2+1$.

Since the CR structure on $\, \Sigma^m$ is integrable, the differential ideal generated by
$\, \{ \, \theta, \, \pi^{0,1} \}$ is closed. By successively differentiating
$\, \pi^{i_{l+1}}_{i_{l-1}} \equiv 0$ $\mod \; \; \theta, \; \pi^{0,1}$ for
$\, l=2, \, ... \, \tau-1$, we obtain the following structure equations.
\begin{equation}\label{fseq}
\pi^{i_{l+1}}_{i_{l}} \wedge \pi^{i_{l}}_{i_{l-1}} \equiv 0  \mod \; \; \theta, \; \pi^{0,1},
\quad   \mbox{for} \; \;  2 \leq \, l \, \leq \, \tau-1.
\end{equation}
Set $\, \pi^{i_{l+1}}_{i_l} \equiv h^{i_{l+1}}_{i_{l} j} \, \pi^j \mod \; \; \theta, \; \pi^{0,1}$.
The sequence of equations (\ref{fseq}) then gives
\begin{equation}
h^{i_{l+1}}_{i_{l} j}h^{i_{l}}_{i_{l-1} k} =  h^{i_{l+1}}_{i_{l} k}h^{i_{l}}_{i_{l-1} j}, \n
\end{equation}
and thus
\begin{equation}
h^{i_{l}}_{i_{l-1} i}h^{i_{l-1}}_{i_{l-2} j} \, ...
\, h^{i_{2}}_{k s} \, = h^{i_{l}}_{ij\, ... \, ks} \n
\end{equation}
is fully symmetric in lower indices.

\begin{definition}
Let  $\, M \hookrightarrow \Sigma^m$ be a CR submanifold of constant
type with the associated decomposition \eqref{decom}
and the Maurer Cartan form \eqref{wstruct}.
The $\, l$-th fundamental form $\, F^l$ is defined by
\begin{align}
F^l &= h^{i_{l}}_{i_{l-1} i}h^{i_{l-1}}_{i_{l-2} j} \, ... \, h^{i_{2}}_{k s} \,
\pi^i \otimes \pi^j \, ... \, \otimes \pi^k \otimes \pi^s \otimes   Z_{i_{l}}  \label{funda} \\
&\equiv \pi^{i_{l}}_{i_{l-1}} \otimes \pi^{i_{l-1}}_{i_{l-2}} \, ... \,
\pi^{i_2}_{i} \otimes \pi^i \otimes  Z_{i_{l}} \mod \theta, \pi^{0,1}.\n
\end{align}
\end{definition}
A  computation with the structure equations (\ref{struct1}) shows
$\, F^l$ is up to scale a well defined section of the bundle $\, S^{l,0}(V^*) \otimes W_l$
over $\, M$, where $\, V$ is the holomorphic or $(1,0)$ tangent space of $\, M$.
$\, F^l$ is one of the simplest $l$-th order extrinsic invariant of the CR submanifold $\, M$.

\subsection{Bochner rigid submanifolds}

The weak structure equations (\ref{wstruct}) is refined in this
section for a class of CR submanifolds whose fundamental forms are
algebraically rigid.

We start with an algebraic definition.
Let $\, V= \mathbb{C}^n, \, W=\mathbb{C}^r$ with the standard Hermitian metric.
Let $\, S^{k,0}(V^*)$ ($S^{0,k}(V^*))$ be the space of holomorphic(anti-holomorphic)
homogeneous polynomials of degree $\, k$,
and denote $\, S^{k,l}= S^{k,0}(V^*) \otimes S^{0,l}(V^*)$.
Let $\, \{ \, x^{i} \, \}$ be a unitary basis of $\, S^{1,0}(V^*)$.
For $\, k, \, l \, \geq 1$, we define the subspace
\begin{equation}
S^{k-1,l-1}_1 = \{ \, x^k \, \bar{x}^k \, f \; \; | \; \; f \in  S^{k-1,l-1} \, \}
                \, \subset S^{k,l}. \n
\end{equation}
Let  $\, \{ \, e_{a} \, \}$ be a unitary basis of $\, W$. For
$\, H=\, H^a \otimes e_a \in S^{k,0}(V^*) \otimes W$
and
$\, P=\, P^a \otimes e_a \in S^{l,0}(V^*) \otimes W$,
we define
\begin{equation}
\langle \, H, \, \bar{P} \, \rangle  = H^a \otimes \bar{P}^a  \in S^{k,l}(V^*). \n
\end{equation}

\begin{definition}\label{boch}
An element $\, H \in S^{k,0}(V^*) \otimes W$ is \emph{Bochner rigid}
if for any $\, P \in S^{k,0}(V^*) \otimes W$, the equation
\begin{equation}
\gamma(H, P) = \langle H, \, \bar{P} \rangle + \langle P, \, \bar{H}\rangle  \in S^{k-1,k-1}_1
\subset S^{k,k} \n
\end{equation}
implies
\begin{equation}\label{boeq}
\gamma(H, P) = 0.
\end{equation}
\end{definition}
Suppose a Bochner rigid $\, H = H^a  \otimes e_a \in S^{k,0}(V^*) \otimes W$ is
nondegenerate, or equivalently the associated map $\, H: \, S^{k,0}(V) \to W$ is surjective.
Then by Cartan's Lemma, (\ref{boeq}) implies
\begin{equation}\label{bochco}
P^a  = u^a_b \, H^b
\end{equation}
for a skew Hermitian matrix $\, u^a_b = - \bar{u}^b_a$.

It is known that when dim $W \leq \frac{1}{2}$ dim $V$, every $\, H
\in S^{k,0}(V^*) \otimes W$ is Bochner rigid [Hu].
We record the following for later application.
\begin{lemma}\label{lemm1}
Let $\, H \in S^{k,0}(V^*) \otimes W$ be a Bochner rigid polynomial.
Suppose for $\, B \in S^{1,0}(V^*) \otimes W$,
\begin{equation}
\langle H, \, \bar{B} \rangle \in S^{k-1,0}_1 \, \subset S^{k,1}.  \n
\end{equation}
Then we necessarily have
\begin{equation}
\langle H, \, \bar{B} \rangle =0. \n
\end{equation}
If $\, H$ is moreover nondegenerate, then  $\, B=0$.
\end{lemma}
\emph{Proof}.$\;$  Let $\, \{ \, x^i \, \}$ be a unitary (1,0)-basis
of $\, V^*$, and denote $\, \langle H, \, \bar{B} \rangle  = (x^i \,
\bar{x}^i) \, g_B$ for $g_B \in S^{k-1,0}$. Take any $\, g \in
S^{k-1,0}$. Then $\, \gamma(H, g\, B) = (x^i \, \bar{x}^i) \,
(g_B \, \bar{g}  + g\, \bar{g}_B ) = 0$ by Bochner rigidity of $\,
H$. Since $\, g$ is arbitrary, $\, g_B=0$.  The rest follows from
the definition of a nondegenerate form.  $\square$

\begin{definition}\label{bordef}
Let $\, M \hookrightarrow \Sigma^m$ be a CR submanifold of constant
type  with the associated decomposition \eqref{decom} and the
fundamental forms \eqref{funda}. $\, M$ is a \emph{Bochner rigid
submanifold} if each of its fundamental forms is Bochner rigid.
\end{definition}

\begin{example}\label{brig}
Let $\, M_{n,p} \hookrightarrow \mathbb{C}P^N = P(\bigwedge^n\mathbb{C}^{n+p})$
be the Pl\"ucker embedding of the Grassmannian
$\, Gr(n, \, \mathbb{C}^{n+p})$, $\, p \leq n$.
Let $\, \hat{M}_{n,p} \hookrightarrow S^{2N+1}$ be the inverse image of
$\, M_{n,p}$ under Hopf map, which is an $\, S^1$-bundle over $\, M_{n,p}$.
When  $\, p = 2, \, 3$,  $\hat{M}_{n,p}$ is a Bochner rigid
CR submanifold of $\, S^{2N+1}$.
\end{example}
We postpone the computation of this example to \textbf{Section 6}.
It is likely the case that $\, \hat{M}_{n,p} \hookrightarrow S^{2N+1}$
is Bochner rigid for all $\, p$.

The main result of this section is the following refined  structure
equations for Bochner rigid submanifolds.
We continue to use the notations of \tb{Section \ref{Fundaforms}}.

\begin{theorem}\label{thm1}
Let $\, M \hookrightarrow \Sigma^m$ be a Bochner rigid submanifold.
Then the weak structure equations \eqref{wstruct} can be refined as
follows.
\begin{align}
\pi_0^A &= 0 \quad   \quad \mbox{for} \; \;  n+1 \leq \, A \, \leq \, m  \label{sstruct} \\
\pi_i^A &= 0 \quad   \quad \mbox{for} \; \;  n_2+1 \leq \, A \, \leq \, m\notag \\
\pi_{i_l}^A &= 0 \quad  \quad \mbox{for} \; \;  n_{l+1}+1 \leq \, A \, \leq \, m,
\quad 2 \leq l \leq \tau-2 \n
\end{align}
These equations furthermore imply
\begin{align}
\pi_{i_l}^{i_{l+1}} &\equiv 0 \quad  \mod \theta, \; \pi^{1,0}
\quad \mbox{for} \; \;  1 \leq \, l \, \leq \, \tau -1, \label{sstruct1} \\
\pi_{m+1}^A &= 0 \quad \; \; \; \mbox{for} \; \;  n_2+1 \leq \, A \, \leq \, m. \n
\end{align}
\end{theorem}
\emph{Proof}.$\;$ By definition of $\, W_2$ and  the weak structure
equations (\ref{wstruct}), we in fact have
\begin{align}
\pi_i^A &= h^A_i \, \theta,  \label{53}  \\
\pi_{m+1}^A &\equiv h^A_i \, \pi^i \quad \mod \theta,
\quad   \quad \mbox{for} \; \;  n_2+1 \leq \, A \, \leq \, m  \n
\end{align}
for a set of coefficients $\, h^A_i$. Set $\, \pi^A_{i_2} \equiv B^A_{i_2j} \bar{\pi}^j$
mod $\theta, \pi^{1,0}$ $\mbox{for} \; \;  n_2+1 \leq \, A \, \leq \, m$.
Differentiating (\ref{53}) mod $\theta$,
\begin{equation}
\im \, h^A_i \, \varpi  \equiv  \pi^A_{i_2} \wedge \pi^{i_2}_{i}
+ \pi^A_{m+1} \wedge \pi^{m+1}_{i} \mod \theta,
\quad \mbox{for} \; \;  n_2+1 \leq \, A \, \leq \, m,\n
\end{equation}
where $\, \varpi= \, \pi^k \w \bar{\pi}^k$, and this gives
\begin{equation}
B^A_{i_2p}\, h^{i_2}_{ij} = -\mbox{i} \,(h^A_i \, \delta_{jp} + h^A_j \, \delta_{ip}).\n
\end{equation}
But $\,h^{i_2}_{ij}$ represents the second fundamental form $\, F^2$ which
is Bochner rigid. By \textbf{Lemma \ref{lemm1}},  $\, B^A_{i_2p} =0,
\, h^A_j=0$ for $n_2+1 \leq \, A \, \leq \, m$. We thus have, since
$\, \pi^A_{i_2} \equiv 0 \mod \theta, \pi^{0,1}$ for $n_3+1 \leq \,
A \, \leq \, m\,$ from (\ref{wstruct}),
\begin{align}
\pi^A_{i} &=0  \quad \quad \quad \quad \quad \quad \quad \,
\quad \quad \mbox{for} \; \;  n_2+1 \leq \, A \, \leq \, m, \n \\
\pi^A_{i_2} &\equiv 0 \quad  \mod \, \theta, \, \pi^{1,0}
\quad \quad \mbox{for} \; \;  n_2+1 \leq \, A \, \leq \, n_3, \n \\
\pi^A_{i_2} &= h^A_{i_2} \, \theta
\quad \quad \quad \quad \quad \quad \quad \; \; \;
\mbox{for} \; \;  n_3+1 \leq \, A \, \leq \, m. \label{57}
\end{align}

Set $\, \pi^A_{i_3} \equiv B^A_{i_3j} \bar{\pi}^j$ mod $\theta, \pi^{1,0}$
for $\, n_3+1 \leq \, A \, \leq \, m$.
Differentiating (\ref{57}) mod $\theta$,
\begin{equation}
h^A_{i_2} \, \varpi  \equiv \pi^A_{i_3} \wedge \pi^{i_3}_{i_2} \mod \theta
\quad \mbox{for} \; \;  n_3+1 \leq \, A \, \leq \, m,\n
\end{equation}
for $\, \pi^A_{m+1} \equiv 0 \mod \theta \, $ for $\, n_2+1 \leq \, A \, \leq \, m$
from \eqref{53}, and this gives
\begin{equation}
B^A_{i_3p}\, h^{i_3}_{i_2 k} = - \mbox{i} \, h^A_{i_2} \, \delta_{kp}.\n
\end{equation}
Multiplying $\, h^{i_2}_{ij}$ and summing over $i_2$ we get
\begin{equation}
B^A_{i_3p}\, h^{i_3}_{i_2 k}\,h^{i_2}_{ij} = - \mbox{i}\, h^A_{i_2}\,h^{i_2}_{ij} \, \delta_{kp}
\quad \mbox{for} \; \;  n_3+1 \leq \, A \, \leq \, m.\n
\end{equation}
But $\,h^{i_3}_{i_2 k}\,h^{i_2}_{ij} = h^{i_3}_{ijk}$ represents  the
third fundamental form $\, F^3$ which is Bochner rigid, and both $\,
F^2, \, F^3$ are nondegenerate by definition. By \textbf{Lemma
\ref{lemm1}}, $\, B^A_{i_3p} =0, \, h^A_{i_2} =0$ for $\, n_3+1 \leq
\, A \, \leq \, m$. We thus have by similar argument as before
\begin{align}
\pi^A_{i_2} &=0   \quad \quad \quad \quad \quad \quad \quad \,
\quad \quad  \mbox{for} \; \; n_3+1 \leq \, A \, \leq \, m, \n \\
\pi^A_{i_3} &\equiv 0 \quad  \mod \, \theta, \, \pi^{1,0}
\quad \quad \mbox{for} \; \;  n_3+1 \leq \, A \, \leq \, n_4, \n \\
\pi^A_{i_3} &= h^A_{i_3} \, \theta
\quad \quad \quad \quad \quad \quad \quad \; \; \;
\mbox{for} \; \;  n_4+1 \leq \, A \, \leq \, m.\n
\end{align}
Continuing in this manner, it is straightforward to verify (\ref{sstruct}).

The first equation of (\ref{sstruct1}) is already proved. Once
(\ref{sstruct}) is true, we have from (\ref{hs}) $\, \pi_{m+1}^A =
p^A \, \theta\,$ $\mbox{for} \; \;  n_2+1 \leq \, A \, \leq \, m\n$.
Differentiating this equation mod $\, \theta$ and collecting
(1,1)-terms,  we get $\, p^A \, \varpi \equiv 0\, $ for   $\, n_2+1
\leq \, A \, \leq \, m.\n \quad \square $

\section{Rigidity of CR submanifolds}

\subsection{CR 1-rigidity}
Let $\, f: \, M \hookrightarrow \Sigma^m$ be a CR submanifold of
constant type. A \emph{type preserving CR deformation} of $\, f$ is
a deformation of $\, f$ through CR immersions with the same type
numbers and the equivalent induced CR structures as $\, f$.
\begin{theorem}\label{thm3}
Let $\, f: \, M \hookrightarrow \Sigma^m$ be a Bochner rigid CR
submanifold  of CR dimension $\geq 2$. Then $\, f$ is 1-rigid under
type preserving CR deformations.
\end{theorem}

The idea is to explore the structure of the compatibility equation
(\ref{ddeforme}) in the presence of the associated decomposition
(\ref{decom}) and the corresponding adapted structure equation
(\ref{sstruct}), (\ref{sstruct1}). The proof is carried out by a
sequence of over-determined pde computations. Simply put, Bochner
rigidity of the fundamental forms implies that each adjacent sequence
of compatibility conditions of the over-determined pde's are tightly
related. It thus plays the role of connecting the propagation of the
associated sequence of compatibility conditions.

It would be interesting and more desirable to remove  the type
preserving hypothesis in the theorem, which is the case when
the height $\,\tau =2$. The computation involved without the hypothesis is, however,
more differential algebraic than algebraic which we are not able to resolve.


\begin{corollary}\label{corn2}{\textnormal{[Hu][EHZ]}}\;
Let $\, f: \, M^n \hookrightarrow \Sigma^m$ be a CR submanifold of constant type
and of CR dimension $\, n \geq 2$. Suppose each of the type
numbers $(\, r_2, \, r_3, \, ... \, r_{\tau} \, )$ are bounded by
$\, r_l \leq \frac{n}{2}$, $\, 2 \leq l \leq \tau$. Then $\, f$ is
1-rigid under type preserving CR deformations.
\end{corollary}

\begin{corollary}\label{cogra}
Let $\, f_{n, p}: \, \hat{M}_{n, p} \hookrightarrow \Sigma^m$ be the
canonical  $\, S^1$-bundle over the Pl\"ucker embedding of the
Grassmannian $\, Gr(n, \, \mathbb{C}^{n+p})$, $p \leq n$.
When $\, p = 2, \, 3$, $\,  f_{n, p}$ is 1-rigid under type preserving CR deformations.
\end{corollary}

We present the proof of the theorem for the case $\, f$ has height $\, \tau =3$,
as the proof for the general case can be easily deduced form this.
We shall agree on the index range
\begin{align}
1 \leq i, \, j, \, k, &\, l, \, p, \,s  \,   \leq n  \n\\
n+1 \leq a, &\, b          \, \leq n_2 \n\\
n_2+1 \leq \nu, &\, \mu    \,\leq n_3, \n
\end{align}
and we continue to use the notations in \textbf{Section 3}.

\emph{Proof of Theorem \ref{thm3}.$\;$} Since $\, M$ is a Bochner
rigid submanifold, the induced Maurer-Cartan  form $\, \pi$ on the
associated $\, H$-bundle $\,B$ takes the following form.
\begin{equation}
\pi =
\begin{pmatrix}
\pi^0_0   & \pi^{0}_{j}    & \pi^{0}_{b}     & 0                & \pi^0_{m+1} \\
\pi^i     & \pi^{i}_{j}    & \pi^{i}_{b}     & 0                & \pi^i_{m+1}\\
0         & \pi^{a}_{j}    & \pi^{a}_{b}     & \pi^{a}_{\mu}    & \pi^a_{m+1} \\
0         & 0              & \pi^{\nu}_{b}   & \pi^{\nu}_{\mu}  & 0 \\
\theta    & \pi^{m+1}_{j}  & 0               & 0                & \pi^{m+1}_{m+1}
\end{pmatrix}
\end{equation}
with
\begin{align}
\pi^{a}_{j}   &= h^a_{jk} \, \pi^k + h^a_j \, \theta, \n \\
\pi^a_{m+1}   &\equiv h^a_j \, \pi^j  \mod \theta, \n \\
\pi^{\nu}_{b} &\equiv h^{\nu}_{bk} \, \pi^k \mod  \theta. \n
\end{align}
Here $\, h^a_{jk}$ and $\, h^{\nu}_{bk} \, h^b_{ij} =
h^{\nu}_{ijk}$  are fully symmetric in lower indices representing
the second and the third fundamental forms. The Lie algebra of $\, H$-valued
group variable $\, V$ in the group action \eqref{gpac} is of the following form
at this stage.
\begin{equation}
V =
\begin{pmatrix}
V^0_0     & V^{0}_{j}    & 0        & 0                &V_{m+1} \\
0         & V^{i}_{j}    & 0        & 0                & V^i_{m+1}\\
0         & 0            & V^{a}_{b}     & 0   & 0 \\
0         & 0            & 0         & V^{\nu}_{\mu}  & 0 \\
0         & 0            & 0        & 0                & V^{m+1}_{m+1}
\end{pmatrix}. \n
\end{equation}

For notational purpose, let $\, \delp$ denote   $\, \pi_1$ in  the
deformation equation (\ref{ddeforme}). Since the type preserving  CR
deformation is considered, it is clear that $\, \delp$ should take
the following form after some absorption by group actions
(\ref{gpac}).
\begin{equation}\label{1delp}
\delp =
\begin{pmatrix}
\delta \pi^0_0   & \delta \pi^{0}_{j}    & \delta\pi^{0}_{b}     & \delta\pi^{0}_{\mu}    & \delta\pi^0_{m+1} \\
0  & \delta\pi^{i}_{j}    & \delta\pi^{i}_{b}     & \delta\pi^{i}_{\mu}                & \delta\pi^i_{m+1}\\
0  & \delta\pi^{a}_{j}    & \delta\pi^{a}_{b}     & \delta\pi^{a}_{\mu}    & \delta\pi^a_{m+1} \\
0  & \delta\pi^{\nu}_{j}  & \delta\pi^{\nu}_{b}   & \delta\pi^{\nu}_{\mu}  & \delta\pi^{\nu}_{m+1} \\
0  & 0                    & 0                     & 0                & \delta\pi^{m+1}_{m+1}
\end{pmatrix}
\end{equation}
Here $\, \delp^i_0$ is translated to $\,0$ by $\, V^i_j-\delta_{ij}V^0_0$ and $\, V^i_{m+1}$
component, and $\, \delp^{m+1}_0$ is translated to $\,0$ by Re$\, V^0_0$ component
for it is a CR deformation, and
$\, \delta\pi^{\nu}_{j} \equiv 0$ mod $\, \theta, \, \pi^{0,1}$ for it is type preserving.
The remaining group variables at this stage are
$\, V^0_0, \, V^a_b, \, V^0_{m+1}$, and  $\, V^{\nu}_{\mu}$
with $\, V^0_0 = - \bar{V}^0_0$.

\emph{Step 1.} Differentiating $\, \delp^i_0 = 0, \, \delp^a_0=0, \, \delp^{\nu}_0 = 0,
\, \delp^{m+1}_0=0$, we get
\begin{align}
\Delta^{i}_{j} \wedge \pi^j   + \,  \delp^i_{m+1} \wedge \, \theta &= 0 \label{72}\\
\delp^{a}_{i} \wedge \pi^i   + \,  \delp^a_{m+1}       \wedge \, \theta &= 0 \n \\
\delp^{\nu}_{i} \wedge \pi^i   + \,  \delp^\nu_{m+1}       \wedge \, \theta &= 0 \n \\
\mbox{Re} \, \delp^0_0 \wedge \, \theta & = 0,\n
\end{align}
where $\, \Delta^i_j = \delp^i_j - \delp^0_0 \, \delta_{ij}$.  Using
the group action (\ref{gpac}) by $\, V^0_{m+1}$ component, we can
absorb $\, \mbox{Re} \, \delp^0_0$ to $\, 0$. Then $\,
(\Delta^{i}_{j})$ is a skew Hermitian matrix valued 1-form such that
\begin{equation}
\Delta^{i}_{j} \wedge \pi^j \equiv 0 \mod \theta,\n
\end{equation}
and Cartan's lemma implies $\, \Delta^{i}_{j} \equiv 0$ mod
$\theta$. $\, \delta\pi^{\nu}_{j} \equiv 0$ mod $\, \theta, \,
\pi^{0,1}$  and (\ref{72}) implies $\, \delta\pi^{\nu}_{j} \equiv 0$
mod $\, \theta$. We thus have
\begin{align}
\begin{pmatrix}
\Delta^{i}_{j} \\
\delp^i_{m+1}
\end{pmatrix}
&=
\begin{pmatrix}
0 & q^i_j  \\
q^i_j & q^i
\end{pmatrix}
\begin{pmatrix}
\pi^j \\
\theta
\end{pmatrix}, \n \\
\begin{pmatrix}
\delp^{a}_{i} \\
\delp^a_{m+1}
\end{pmatrix}
&=
\begin{pmatrix}
p^a_{ij} & p^a_i \\
p^a_i & p^a
\end{pmatrix}
\begin{pmatrix}
\pi^j \\
\theta
\end{pmatrix}, \n\\
\begin{pmatrix}
\delp^{\nu}_{i} \\
\delp^{\nu}_{m+1}
\end{pmatrix}
&=
\begin{pmatrix}
0 & p^{\nu}_i \\
p^{\nu}_i & p^{\nu}
\end{pmatrix}
\begin{pmatrix}
\pi^j \\
\theta
\end{pmatrix}, \n \\
\mbox{Re} \, \delp^0_0  &= 0,\n
\end{align}
where $\, q^j_i = - \bar{q}^i_j$, $\, p^a_{ij} = p^a_{ji}$.
Differentiating $\, \mbox{Re} \, \delp^0_0 = 0$ with these
relations, we get
\begin{equation}
\delp^0_{m+1} = \frac{\mbox{i}}{2} ( \, \bar{q}^i \, \pi^i - \, q^i \, \bar{\pi}^i \, )
+ \, q^0 \, \theta,\n
\end{equation}
for a real coefficient $\, q^0$.
Note $\, \delta \pi^{m+1}_{m+1}  = - \delta \bar{\pi}^{0}_{0} = \delta \pi^{0}_{0}$.
The remaining group variables at this stage are
$\, V^0_0, \, V^a_b$, and  $\, V^{\nu}_{\mu}$.

\emph{Step 2.}  Differentiating $\, \Delta^i_j = q^i_j \, \theta$
mod $\, \theta$ and collecting terms, we obtain an equation which is
equivalent to
\begin{equation}\label{82}
h^a_{jk} \, \bar{p}^a_{il} \, + \, p^a_{jk} \,  \bar{h}^a_{il}
= \, \mbox{i} \, ( q^i_j \, \delta_{kl} + \, q^i_k \, \delta_{jl}
+  \, q^l_j \, \delta_{ik} + \, q^l_k \, \delta_{ij}).
\end{equation}
Since $\, h^a_{ij}$ is Bochner rigid, this equation implies
\begin{align}
q^i_j    &= 0, \n \\
h^a_{jk} \, \bar{p}^a_{il} \, + \, p^a_{jk} \,  \bar{h}^a_{il}   &= 0. \n
\end{align}
From (\ref{bochco}), $\, p^a_{ij} = \, u^a_b \, h^b_{ij}$ for  a skew
Hermitian matrix $\, (u^a_b)$, and thus can be absorbed to $\, 0$ by
the group action by $\, V^a_b-\delta_{ab} V^0_0$ component. Differentiating
$\, \delp^i_{m+1} = q^i \, \theta$ mod $\theta$,
\begin{align}
q^i \, \varpi &\equiv \pi^i \wedge \delp^0_{m+1} + \,\pi^i_a \wedge \delp^a_{m+1} \quad \mod \theta \n\\
&\equiv \pi^i  \wedge \,
( \, \frac{\mbox{i}}{2} ( \, \bar{q}^j \, \pi^j - \, q^j \, \bar{\pi}^j \, ))
- \bar{h}^a_{ik} \, p^a_j \, \bar{\pi}^k \wedge \pi^j \quad \mod \theta,\n
\end{align}
where $\, \varpi= \pi^k \wedge \bar{\pi}^k$.
Collecting $(2,0)$-terms, we get $\, q^i = 0$. Since $\, h^a_{ij}$
represents the nondegenerate second fundamental form, the remaining
equation $\, \bar{h}^a_{ik} \, p^a_j = 0$ implies $\, p^a_i$ = 0.
Now differentiating $\, \pi^a_{m+1} = p^a \, \theta$ mod $\theta$,
we get $\, p^a =0$. $\, \delp$ is now reduced to the
following form.
\begin{equation}
\delp =
\begin{pmatrix}
\delta \pi^0_0   & 0   & 0    & \delta\pi^{0}_{\mu}    & q^0 \, \theta \\
0  & \delp^0_0 \, \delta_{ij}   & 0     & \delta\pi^{i}_{\mu}                & 0\\
0  & 0    & \delta\pi^{a}_{b}     & \delta\pi^{a}_{\mu}    & 0 \\
0  & \delta\pi^{\nu}_{j}  & \delta\pi^{\nu}_{b}   & \delta\pi^{\nu}_{\mu}  & \delta\pi^{\nu}_{m+1} \\
0  & 0                    & 0                     & 0                & \delta \pi^0_0
\end{pmatrix}\n
\end{equation}
The remaining group variables at this stage are
$\, V^0_0$ and  $\, V^{\nu}_{\mu}$.

\emph{Step 3.}  Differentiating $\, \delp^{\nu}_{i}= p^{\nu}_i \,
\theta$  mod $\theta$, we get
\begin{equation}\label{88}
p^{\nu}_i \, \varpi \equiv \delp^{\nu}_a \wedge \, \pi^a_i \quad \mod \theta.
\end{equation}
Set $\, \delp^{\nu}_a \equiv p^{\nu}_{ai} \, \pi^i + C^{\nu}_{ai} \, \bar{\pi}^i
\mod \theta$.
Collecting $(2,0)$, and $(1,1)$ terms in (\ref{88}) we get
\begin{align}
p^{\nu}_{ai} \, h^a_{jk} &= p^{\nu}_{aj} \, h^a_{ik} = p^{\nu}_{ijk}
\quad \mbox{fully symmetric in lower indices} \n\\
\mbox{i} \, p^{\nu}_i \, \delta_{jk} &= - C^{\nu}_{ak} \, h^a_{ij}. \label{322}
\end{align}
By Bochner rigidity of $\, h^a_{ij}$, (\ref{322}) implies  $\,
C^{\nu}_{ak} = 0$, $\, p^{\nu}_i = 0$. Differentiating $\,
\delp^{\nu}_{m+1} = p^{\nu} \, \theta \,$ mod $\, \theta$ with these relations
gives $\, p^{\nu} =0$. Differentiating $\, \delp^0_{m+1} = \, q^0 \, \theta$
mod $\theta$ at this stage gives $\, q^0 = 0$.
$\, \delp$ is now  reduced to the following form.
\begin{equation}
\delp =
\begin{pmatrix}
\delta \pi^0_0   & 0   & 0    & 0    & 0 \\
0  & \delp^0_0 \, \delta_{ij}   & 0     & 0               & 0\\
0  & 0    & \delta\pi^{a}_{b}     & \delta\pi^{a}_{\mu}    & 0 \\
0  &0  & \delta\pi^{\nu}_{b}   & \delta\pi^{\nu}_{\mu}  & 0 \\
0  & 0                    & 0                     & 0                & \delta \pi^0_0
\end{pmatrix}\n
\end{equation}

\emph{Step 4.}  Differentiating $\, \delp^{a}_{i}=0$ mod $\theta$,
\begin{equation}
\Delta^a_b \, \wedge \pi^b_i \equiv 0 \quad \mod \theta \n
\end{equation}
where $\, \Delta^a_b = \delp^a_b - \delp^0_0 \, \delta_{ab}$. Since
the coefficients of   $\, \pi^a_{i} \equiv h^a_{ij} \, \pi^j \; \mod
\theta \;$ represent the second fundamental form which is
nondegenerate by definition, a variant of Cartan's lemma, [GH],
implies $\, \Delta^{a}_{b} \mod \theta$ is in the span of $\, \{ \,
\pi^a_i\}$. But $\, \pi^a_i$ consists of $(1,0)$-forms mod $\,
\theta$ and $\, \Delta^{a}_{b}$ is skew Hermitian. Hence $\,
\Delta^{a}_{b} \equiv 0$ mod $\theta$, and we write
\begin{equation}
\Delta^a_b = q^a_b \, \theta \n
\end{equation}
where $\, q^b_a = -\bar{q}^a_b$. Differentiating this equation,
we get $\, - \, \im \,q^a_b \, \varpi  \equiv \bar{\pi}^{\nu}_a
\wedge \delp^{\nu}_b + \delta \bar{\pi}^{\nu}_a \wedge \pi^{\nu}_b
\mod \theta $, or equivalently
\begin{equation}
- \,\mbox{i} \, q^a_b \, \delta_{sk}
= \bar{h}^{\nu}_{as} \, p^{\nu}_{bk} +  \bar{p}^{\nu}_{as} \, h^{\nu}_{bk}. \n
\end{equation}
Multiplying both sides by $\, \bar{h}^a_{lp} \, h^b_{ij}$ and
summing  with respect to $\, a, \, b$,
\begin{equation}
- \mbox{i} \, q^a_b \, \bar{h}^a_{lp} \, h^b_{ij} \delta_{sk}
= \bar{h}^{\nu}_{lps} \, p^{\nu}_{ijk} +  \bar{p}^{\nu}_{lps} \, h^{\nu}_{ijk}. \n
\end{equation}
By Bochner rigidity of $\, h^{\nu}_{ijk}$ representing  the third
fundamental form, this implies
\begin{align}
q^a_b &=0,\n \\
p^{\nu}_{ijk} &= p^{\nu}_{ai}\, h^a_{jk} = u^{\nu}_{\mu} \, h^{\mu}_{ijk}  \n \\
&= u^{\nu}_{\mu} \, h^{\mu}_{ai}\, h^a_{jk} \n
\end{align}
for a skew Hermitian $\, u^{\mu}_{\nu}= - \bar{u}^{\nu}_{\mu}$.
Since $\, h^a_{ij}$ is nondegenerate, this implies
\begin{equation}
p^{\nu}_{ai} = u^{\nu}_{\mu} \, h^{\mu}_{ai}, \n
\end{equation}
and we may absorb $\, p^{\nu}_{ai}$ to $\, 0$ by group action by
$\, V^{\nu}_{\mu}-\delta_{\nu \mu} V^0_0$ component as before.
Put $\, \delp^{\nu}_{a} = p^{\nu}_{a} \, \theta$. Differentiating
$\, \Delta^a_b=0$ and collecting $\, \theta \, \wedge \, \pi^{1,0}$-terms,
we get $\, h^{\nu}_{ai} \, \bar{p}^{\nu}_b = 0$. Multiplying $\, h^a_{jk}$ and
summing with respect to $\, a$, we get
\begin{equation}
h^{\nu}_{ijk} \, \bar{p}^{\nu}_b = 0. \n
\end{equation}
Since $\, h^{\nu}_{ijk}$ represents the third fundamental form
which is nondegenerate by definition, this equation implies $\,
p^{\nu}_b = 0$.  $\, \delp$ is now reduced to the following form.
\begin{equation}
\delp =
\begin{pmatrix}
\delta \pi^0_0   & 0   & 0    & 0    & 0  \\
0  & \delta \pi^0_0 \, \delta_{ij}   & 0     & 0               & 0   \\
0  & 0    & \delta \pi^0_0 \, \delta_{ab}     & 0    & 0  \\
0  &0  & 0   & \delta\pi^{\nu}_{\mu}  & 0  \\
0  & 0                    & 0                     & 0                & \delta \pi^0_0
\end{pmatrix}.\n
\end{equation}
Note there are no remaining group variables at this stage.

\emph{Step 5.}  Differentiating $\, \delp^{\nu}_{a}=0$ mod $\theta$,
\begin{equation}
\Delta^{\nu}_{\mu} \, \wedge \pi^{\mu}_a \equiv 0 \quad \mod \theta \n
\end{equation}
where $\, \Delta^{\nu}_{\mu} = \delp^{\nu}_{\mu} - \delp^0_0 \, \delta_{\nu \mu}$.
Since the coefficients of  $\, \pi^{\mu}_a \equiv h^{\mu}_{aj} \, \pi^j$ mod $\, \theta$
represents the nondegenerate third fundamental form,
a variant of Cartan's lemma as before implies
\begin{equation}
\Delta^{\nu}_{\mu} = \delp^{\nu}_{\mu} - \delp^0_0 \, \delta_{\nu \mu}
 \equiv 0 \quad \mod \theta. \n
\end{equation}

Now $\, \delp$ takes values in $\, \mathfrak{su}(m+1,1)$,  and in
particular tr$\, \delp =0$. Hence $\, \delp^0_0$ is a multiple of
$\, \theta$, and consequently
\begin{equation}
\delp = \mbox{P} \, \theta \n
\end{equation}
for a scalar $\, \mathfrak{su}(m+1,1)$-valued coefficient $\, \mbox{P}$.
But differentiating this equation mod $\,\theta$ gives
$\, \mbox{P} \, \varpi \equiv 0$, hence $\, \mbox{P} =0$. $\; \square$

Here we summarize the process of computation for 1-rigidity in
analogy with the game of dominoes. Imagine $\, \delp$ is a
$\, (\tau + 2)$-by-$\,(\tau + 2)$ standing blocks of dominoes. Our objective is to
lay down all of these domino blocks. There are three main ingredients in
the process.

\tb{A}. The initial push: This corresponds to imposing  the
condition of admissible deformations as in (\ref{1delp}).

\tb{B}. When a block gets hit, it falls toward the nearby standing
blocks: This is expressed by the compatibility conditions obtained
by differentiation using the deformation equation (\ref{ddeforme}).

\tb{C}. The nearby standing blocks in process \tb{B} are within
hitting distance, and gets hit by at least one domino block, and
falls: This in CR case is guaranteed by the Bochner rigid
conditions.

Once there was the initial push \tb{A}, the processes \tb{B}  and
\tb{C} alternate just as it is the case with real dominoes. We
finally mention that in process \tb{C}, it is conceivable that a
standing block gets hit but the force is not enough to lay it down,
or it is equally conceivable that a standing block gets hit by more
than one blocks and the hitting forces cancel. In these cases, one
has to analyze the prolongation of the deformation equation. This
phenomena occurs when trying to remove the type preserving
hypothesis from \tb{Theorem \ref{thm3}}.

\subsection{Local rigidity}
Although it is formulated based on a natural geometric consideration,
1-rigidity introduced in \textbf{Section 2} is an abstract
definition. It is not clear whether 1-rigidity has any bearing on
the actual rigidity property of submanifolds under admissible
deformation. The main result of this section is that 1-rigidity can
be extended to actual local rigidity at least for Bochner rigid CR
submanifolds in spheres. The idea of proof is a certain
averaging trick, which might be of use in rigidity problems in other
homogeneous spaces as well.

\begin{theorem}\label{lorig}
Let $\, f_0: \, M \hookrightarrow \Sigma^m$ be a Bochner rigid
CR submanifold  of height $\tau$. Then there exists a
$\, C^{\tau}$ neighborhood $\, \mathcal{U} \subset  C^{\tau}(M,
\, \Sigma^m)$ of $\, f_0$ with the following property. Suppose $\, f,
\, g \in \mathcal{U}$ are CR immersions such that their induced CR
structures on $\, M$ are mutually equivalent and that they are of
the same type as $\, f_0$. Then $\, f$ and  $\, g$ are congruent up to an
automorphism of $\, \Sigma^m$.
\end{theorem}
Note when $\, \tau=2$, the type preserving hypothesis can be dropped.

We present the proof for the case $\, \tau = 2$, for the proof for
the general case is essentially the same. We follow the notations
in \tb{Section 4.1}.

\emph{Proof of Theorem \ref{lorig}.$\;$} Given two CR immersions $\, f,
\, g$,  let $\, F: \, B_f \hookrightarrow f^*$SU$(m+1,1)$, $\, G: \,
B_g \hookrightarrow  g^*$SU$(m+1,1)$ be the associated adapted $\,
H$-bundles over $\, M$. Let
\begin{align}
s_f &= ( F_0, \, F_1, \, ... \, F_{m+1} ), \n \\
s_g &= ( G_0, \, G_1, \, ... \, G_{m+1} )\n
\end{align}
be any sections of these bundles, and denote
\begin{align}
s_f^*(F^*\phi) &= \alpha,  \n \\
s_g^*(G^*\phi) &= \beta \n
\end{align}
where $\, \phi$ is the Maurer-Cartan form of SU$(m+1,1)$.
Set
\begin{align}
\pi &= \frac{1}{2}( \alpha + \beta ), \notag \\
\delta \pi &= \alpha - \beta, \n
\end{align}
and observe
\begin{align}
- \, d (\delta \pi) &=  \alpha \wedge \alpha  - \beta  \wedge \beta   \n  \\
&=  \delta \pi \wedge \pi + \pi \wedge \delta \pi. \label{111}
\end{align}
For a fixed section $\, s_g$, we wish to show there exists a section
$\, s_f$ such that $\, \delp = 0$.

Since $\, f$ and $\, g$ induce the equivalent CR structures on $\, M$, we may
take a section $\, s_f$ such that
\begin{align}
\alpha_0^i &= \beta_0^i = \pi^i,   \label{113}\\
\alpha_0^{m+1} &= \beta_0^{m+1} = \theta \n
\end{align}
and $\, \delta \pi$ takes the following form.
\begin{equation}
\delta \pi =
\begin{pmatrix}
\delta \pi^0_0 & * & *  & \delta \pi^0_{m+1} \\
0 & * & *  & *\\
0 & * & *  & * \\
0 & 0 & 0  & *
\end{pmatrix}. \n
\end{equation}
Differentiating $\, \delta \pi^{m+1}_0 =0$ using (\ref{111}), we get
Re$\, \delta \pi^0_0 \wedge \theta =0$. Modifying $\, F_{m+1}$ by a
multiple of $\, F_0$ if necessary, which is still a legitimate
section of $\, B_f$ preserving the relations (\ref{113}), we may in
fact have Re$\, \delta \pi^0_0 =0$(this requires a computation,
which is rather similar to the absorption by group action as in
CR 1-rigidity and we omit). Now differentiating $\, \delta \pi^i_0 =0$
we get
\begin{equation}\label{115}
\begin{pmatrix}
\Delta^i_j \\
\delta \pi^i_{m+1}
\end{pmatrix}
=
\begin{pmatrix}
0 & q^i_j \\
q^i_j & q^i
\end{pmatrix}
\begin{pmatrix}
\pi^j \\
\theta
\end{pmatrix}
\end{equation}
where $\, \Delta^i_j = \delp^i_j - \delp^0_0 \, \delta_{ij}$,  $\,
q^j_i = - \, \bar{q}^i_j$.

Let $\, f^a = f^a_{ij}\, \pi^i \, \circ \pi^j$ and  $\, g^a =
g^a_{ij}\, \pi^i \,\circ \pi^j$ represent the second fundamental
forms of $\, f$ and  $\, g$. If $\, f$ and $\, g$ are sufficiently
$\, C^2$ close to $\, f_0$  which is Bochner rigid, then $\, f$ and
$\, g$ would be Bochner rigid. Thus the average $\, \frac{1}{2} \, (
f^a + g^a )$ would be Bochner rigid too for a section $\, s_f$ with
appropriate $\, (F_{n+1}, \, F_{n+2}, \, ... \, F_{n+r})$-part.
Differentiating (\ref{115}) and proceeding as in  (\ref{82}),  we
get
\begin{align}
q^i_j &=0, \n \\
f^a - g^a &= u^a_b \, (  f^a + g^a ) \label{116}
\end{align}
for a skew Hermitian matrix $\, (u^a_b) = u =  -u^*$.

Set $\vec{f} = (\, f^{n+1}, \, f^{n+2}, \, ... \, f^{n+r} \, )^t$, $\vec{g} =
(\, g^{n+1}, \, g^{n+2}, \, ... \, g^{n+r} \, )^t$,
and  write (\ref{116}) in matrix form
\begin{equation}
(I - u ) \, \vec{f} = (I + u ) \, \vec{g}. \n
\end{equation}
When $\, f$ and $\, g$ are sufficiently $\, C^2$ close to  $\, f_0$,
$\, u$ is small and $\, I-u$ would be invertible. Thus we may write
\begin{equation}
\vec{f} = (I - u)^{-1} \, (I + u ) \, \vec{g}. \n
\end{equation}
But
\begin{align}
(I - u)^{-1} \, (I + u )((I - u)^{-1} \, (I + u ))^*
&= (I - u)^{-1} \, (I + u )((I - u) \, (I + u )^{-1}  \n  \\
&= (I - u)^{-1} \, (I - u )((I + u) \, (I + u )^{-1} \n   \\
&= I. \n
\end{align}
Thus we may modify $\, (F_{n+1}, \, F_{n+2}, \, ... \, F_{n+r})$
part of the section $\, s_f$ by unitary matrices close to identity
to have
\begin{equation}
\delta \pi^a_i = 0. \n
\end{equation}
$\, \delta \pi$ now becomes
\begin{equation}
\delta \pi =
\begin{pmatrix}
\delta \pi^0_0 & * & *  & \delta \pi^0_{m+1} \\
0 & \delp^0_0 \, \delta_{ij} & 0  & * \\
0 & 0 & *  & * \\
0 & 0 & 0  & \delta \pi^0_0
\end{pmatrix}. \n
\end{equation}

The rest of the computation proceeds the same as in the proof of
\textbf{Theorem \ref{thm3}}, and we conclude for this section $\, s_f$,
\begin{equation}
\delta \pi = \alpha - \beta = 0. \n
\end{equation}
By the fundamental theorem [Gr, p\,780], the sections $\, s_f$  and
$\, s_g$ are congruent by an element of SU($m+1,1$),
and hence $\, f$ and $\, g$ are congruent by an automorphism of $\, \Sigma^m$.
$\; \square$

\section{Whitney submanifold}

The main result of this section is the local characterization that
every nonlinear CR-flat submanifold $\, M^n \subset \Sigma^{2n}$, $\, n \geq 2$,
is a part of a \emph{Whitney submanifold}.
Whitney submanifold is an example of a CR submanifold
which is not 1-rigid. The result of this section in fact implies
it is CR deformable in exactly 1 direction, see the end of this section.

Let $\, V^{n+2}= \mathbb{C}^{n+1,1}$ be the complex vector space
with coordinates $\, \xi = (\, \xi^0, \, \xi^i, \, \xi^{n+1} \,)$,
$\, 1 \, \leq i \, \leq  n$, and a Hermitian scalar product
\begin{equation}
Q_n( \, \xi, \, \bar{\xi} \, ) = \xi^i \, \bar{\xi}^i +
\mbox{i} \, (\xi^0 \, \bar{\xi}^{n+1} - \xi^{n+1} \, \bar{\xi}^0 ).\n
\end{equation}
Let $\, \Sigma^n \simeq S^{2n+1}$ be the set of equivalence classes up to scale of
null vectors.

Let $\, \mu = (\, \mu^0, \, \mu^i, \, \mu^{n+i} \,
\mu^{2n+1} \, )$, $1 \leq i \leq n$, be  the coordinates of  $\, V^{2n+2}$.
\emph{Whitney submanifold} $\, \Gamma_n: \, \Sigma^n \to
\Sigma^{2n}$  is an immersion induced by the quadratic map $\,
\widehat{\Gamma}_n: \, V^{n+2} \to  V^{2n+2}$ defined as
\begin{align}
\mu^0 &= 2 \, \xi^0 \, \xi^{n+1}, \n \\
\mu^{2n+1} &= \, (\xi^{n+1})^2 - (\xi^{0})^2, \n \\
\mu^i &= \xi^i  \, (\, \mbox{i} \, \xi^0 + \xi^{n+1}), \n \\
\mu^{n+i} &= \xi^i \, (-\, \mbox{i} \, \xi^0 + \xi^{n+1}). \n
\end{align}
$\,  \widehat{\Gamma}_n^* \, Q_{2n} = 2 \, (\, \xi^0 \, \bar{\xi}^0
+ \xi^{n+1} \, \bar{\xi}^{n+1}) \, Q_n$,  and the induced map
$\, \Gamma_n$ on $\, \Sigma^n$ is well defined.
It is easy to check $\, \Gamma_n$ is
CR-equivalent to the boundary map $\, \partial W_n: \, S^{2n+1} \, \to \, S^{4n+1}$
of the following \emph{Whitney map}
$\, W_n: \, \mathbb{B}^{n+1} \, \to \, \mathbb{B}^{2n+1}$, where
$ \,  \mathbb{B}^{\star} \subset \mathbb{C}^{\star}$ is the unit ball
and $\, (\, z^0, \, z^i \, )$, $\, 1 \, \leq i \, \leq  n$, is
a coordinate of $\, \mathbb{C}^{n+1}$.
\begin{equation}\label{whitneymap}
W_n(z^0, z^i) = (\, (z^0)^2, \, z^0 \, z^i, \, z^i \, ).
\end{equation}
This equivalence is via the isomorphism $\, \Sigma^n \simeq
S^{2n+1}$  given in coordinates
\begin{align}
z^0 &= \frac{ \mbox{i} \, \xi^0 + \xi^{n+1} }{ -\, \mbox{i} \, \xi^0 + \xi^{n+1}}  \n \\
z^i &= \frac{\sqrt{2} \, \xi^i}{-\, \mbox{i} \, \xi^0 + \xi^{n+1}}  \n
\end{align}
Set $\, \Sigma^n_0 = \{ \, [\xi] \in \Sigma^n \, | \;    \xi^i = 0, \; \forall \, i \,\}$
and $\, \Sigma^n_s = \{ \, [\xi] \in \Sigma^n \, | \;  \mbox{i} \,
\xi^0 + \xi^{n+1} = 0 \, \}$. Then $\, \Gamma_n$ is an
immersion which is 1 to 1 on $\, \Sigma^n - \Sigma^n_0$, 2 to 1 on
$\, \Sigma^n_0 $, and  the second fundamental form vanishes
along $\, \Sigma^n_s$.

\begin{theorem}\label{whit}
Let $\, M^n \hook \Sigma^{2n}$ be a $\, C^3$ nonlinear CR-flat
submanifold of CR dimension and codimension $\, n \geq 2$.
Then $\, M$ is congruent to a part of the Whitney
submanifold up to an automorphism of $\, \Sigma^{2n}$.
\end{theorem}
CR-flat submanifold $\, M^1 \hook \Sigma^2$ has been classified
by Faran [Fa]. In contrast to $\, n \geq 2$ cases, there are four
inequivalent CR flat submanifolds when $\, n = 1$.

As a corollary, we have a simple characterization of the proper holomorphic
maps from a unit ball $\, \mathbb{B}^{n+1}$ to $\, \mathbb{B}^{2n+1}$ [HJ].
\begin{corollary}\label{whitco}
Let $\, F: \, \mathbb{B}^{n+1} \, \to \, \mathbb{B}^{2n+1}$, $n \geq 2$,
be a nonlinear proper holomorphic map which is $\, C^{n+1}$ up to the
boundary. Then $\, F$ is equivalent to the Whitney map \eqref{whitneymap}
up to automorphisms of the unit balls.
\end{corollary}

\textbf{Theorem \ref{whit}} is based on the following algebraic
lemma due to Iwatani on the normal form of the second fundamental form
of a Bochner-K\"ahler submanifold [Iw][Br3]. Let $\, V = \mathbb{C}^n$,
$\, W = \mathbb{C}^n$ with the standard Hermitian scalar product.
Let $\,\{ \,  z^i \, \}$ be a unitary (1,0)-basis for $\, V^*$, and
let $\,\{ \,  w_a \, \}$ be a unitary basis for $\, W$.
\begin{lemma}\label{iwa}{\textnormal{[Iw]}}\;
Suppose $\, H=\, h^a_{ij} \, z^i \, z^j \otimes w_a \in S^{2,0}(V^*)\otimes W$ satisfies
\begin{equation}
\gamma(H, H) = \, h^a_{ij} \, \bar{h}^a_{kl}
\, z^i \, z^j \otimes \bar{z}^k \, \bar{z}^l
\in S^{1,1}_1 \subset S^{2,2},\n
\end{equation}
or simply $\, \gamma(H, H)$ is Bochner-flat [Br3].
Then up to a unitary transformation on $\, V$,
\begin{equation}
H=\, h^a_{in} \, z^i \, z^n \otimes w_a.\n
\end{equation}
\end{lemma}

Set $\, \nu_i = h^a_{in} \, w_a \in W$. A computation shows  $\,
\gamma(H, H)$ is Bochner-flat whenever
\begin{align}
\langle\, \nu_i,  \, \nu_j \, \rangle &= 0 \quad \mbox{for} \; i \ne j, \notag \\
\langle\, \nu_n,  \, \nu_n \, \rangle &= 4 \, <\, \nu_q,  \, \nu_q \, > \quad
\mbox{for} \; \, q < n.\n
\end{align}
Thus up to a unitary transformation on $\, W$, we may assume
\begin{align}
\nu_q &= r \, w_q \quad \mbox{for} \, \; q < n, \notag \\
\nu_n &= 2\,r \, w_n, \n
\end{align}
for some $\, r \geq 0$.

Let $\, M^n \hookrightarrow \Sigma^{2n}$ be a CR-flat submanifold.
The second fundamental form of $\, M$ is Bochner-flat [EHZ].
In the notation of \textbf{Section 3}, \textbf{Lemma \ref{iwa}}
then implies we may write
\begin{align}
\pi^{n+i}_q &\equiv r \, \delta_{iq} \, \omega^n \quad \mod \; \theta,
\quad \mbox{for} \; q < n, \notag \\
\pi^{n+i}_n &\equiv r (1+\delta_{in}) \, \omega^i \quad \mod \; \theta, \n
\end{align}
where $\, \theta$ is the dual to the contact hyperplane fields.
Assume $\, M$ is not linear,
$\, H \ne 0$, and we scale $\, r=1$ using the group action by Re$\, \pi_0^0$.
We thus obtain the following local normal form of the second fundamental
form of a nonlinear CR-flat submanifold $\, M^n \hookrightarrow \Sigma^{2n}$.
\begin{align}
\pi^{n+i}_q &= \delta_{iq} \, \omega^n + h^i_q \theta
\quad \mbox{for} \; q < n,  \label{rs1} \\
\pi^{n+i}_n &= (1+\delta_{in}) \, \omega^i + h^i_n \theta \n
\end{align}
for some coefficients $\, h^i_j$.

\textbf{Theorem \ref{whit}} is now obtained by successive
differentiation of this normalized structure equation. We assume
$\, n \geq 3$ for simplicity for the rest of this section,  as
$\, n=2$ case can be treated in a similar way. We shall agree on the
index range $\, 1 \leq p, \, q, \, s, \, t \, \leq n-1$,  and denote
$\, p'=n+p$, $\, n'=n+n$. Recall our convention $ - \, d \theta \equiv \, \mbox{i} \,
\pi^k \w \, \bar{\pi}^k \n \equiv  \, \mbox{i} \, \varpi \mod \theta $,
and we denote $\, \pi^i = \omega^i$ for the sake of notation.

\emph{Step 1. $\;$} Differentiating $\, \pi^{n+n}_s  =  h^n_s \, \theta$ $\mod{\theta}$,
we get
\begin{equation}
\mbox{i}\,h^n_s \, \varpi \equiv ( \pi^{n'}_{s'} - 2 \pi^n_s) \wedge \omega^n
+ \pi^{n'}_{2n+1} \wedge (-\, \im \, \bar{\omega}^s) \quad \mod \theta.  \n
\end{equation}
Since $\, n-1 \geq 2$, this implies $\, h^n_s = 0$, and by Cartan's lemma
\begin{equation}
\begin{pmatrix}
\pi^{n'}_{s'} - 2 \pi^n_s \\
\pi^{n'}_{2n+1}
\end{pmatrix}
\equiv
\begin{pmatrix}
c_s & u \\
u  & 0
\end{pmatrix}
\begin{pmatrix}
\omega^n \\
-\im \, \bar{\omega}^s
\end{pmatrix}  \quad \mod{\theta}\n
\end{equation}
for coefficients $\, c_s, \, u$.

\emph{Step 2. $\;$}  Differentiating  $\, \pi^{t'}_s = h^t_s \,
\theta$ $\mod{\theta}\, $ for $\, t \ne s$, we get
\begin{equation}
\mbox{i}\,h^t_s \, \varpi \equiv ( \pi^{t'}_{s'} -  \pi^t_s) \wedge \omega^n - \pi^n_s \wedge \omega^t
+ \pi^{t'}_{2n+1} \wedge (-\im \, \bar{\omega}^s) \quad \mod \theta.\n
\end{equation}
Since $\, n-1 \geq 2$, this implies $\, h^t_s = 0$ for $\, t \ne s$,
and by Cartan's lemma
\begin{equation}
\begin{pmatrix}
\pi^{t'}_{s'} -  \pi^t_s \\
-\pi^n_s \\
\pi^{t'}_{2n+1}
\end{pmatrix}
\equiv
\begin{pmatrix}
0 & b_s   & -\im \, \bar{b}_t \\
b_s  & 0  & a \\
-\im \, \bar{b}_t  & a  & 0
\end{pmatrix}
\begin{pmatrix}
\omega^n \\
\omega^t \\
-\im \, \bar{\omega}^s
\end{pmatrix} \quad \mod{\theta}\n
\end{equation}
for coefficients $\, c_s, \, u$.
Since $\, \pi^{t'}_{s'} -  \pi^t_s$ is skew Hermitian, it cannot have
any $\, \omega^n$-term.

\emph{Step 3. $\;$}  Differentiating  $\, \pi^{t'}_t = \omega^n +
\mbox{i}\,h^t_t \, \theta \,$ $\mod{\theta}$ and collecting terms, we get
\begin{equation}
h^t_t \, \varpi \equiv ( \pi^{t'}_{t'} -  \pi^t_t+ \pi^0_0 - \pi^n_n) \wedge \omega^n
+( b_p \, \omega^n - \im \, a  \, \bar{\omega}^p ) \wedge \omega^p
+( -b_t \, \omega^t + \bar{b}_t \, \bar{\omega}^t ) \, \wedge \omega^n \n
\end{equation}
$\mod \theta$. Since $\, n-1 \geq 2$, this implies $\, h^t_t = a$, and
\begin{align}
\Delta_t   &= \pi^{t'}_{t'} -  \pi^t_t+ \pi^0_0 - \pi^n_n \notag \\
&=  a_t \, \omega^n - \im \, a \, \bar{\omega}^n + ( b_t \, \omega^t - \bar{b}_t \, \bar{\omega}^t )
+ \sum_p \, b_p \,  \omega^p  - A_t \, \theta \n
\end{align}
for coefficients $\, a_t, \, A_t$.

\emph{Step 4. $\;$}  From \emph{Step 2},  we may use the group
action by $\, \pi_{2n+1}^n$ to translate $\, a = 0$, which we assume
from no on. We also translate $\, h^t_n = 0$ by $\, \pi_{2n+1}^t$.
Differentiating $\, \pi^{t'}_n = \omega^t$ $\mod{\theta}\, $ with
these relations, we get
\begin{equation}
0 \equiv \bar{b}_t \, ( \, \sum_p \omega^p \wedge \bar{\omega}^p )
+( 2 \bar{c}_t - 3 \bar{b}_t \,) \omega^n \wedge \bar{\omega}^n
+ \omega^n \wedge ( a_t + 2 \im \,\bar{u} ) \wedge \omega^t \quad \mod{\theta}. \n
\end{equation}
Thus $\, b_t = c_t = 0$, $\, a_t = - 2 \, \im \,\bar{u}$.

\emph{Step 5. $\;$}   Differentiating $\, \pi^{n'}_n = 2 \, \omega^n
+ h^n_n \, \theta$ $\mod{\theta}$ and collecting terms, we get $\,
A_t = \, A$ for a single variable, and
\begin{equation}
\mbox{i}\, h^n_n \, \varpi \equiv 2 \, ( \pi^{n'}_{n'} -  \pi^n_n+ \pi^0_0 - \pi^n_n) \wedge \omega^n
+  \, \im \, u \, \omega^p \wedge \bar{\omega}^p
- \, \im \, u \, \omega^n \wedge \bar{\omega}^n \quad \mod{\theta}.\n
\end{equation}
This implies $\, h^n_n = \, u$, and
\begin{align}
\Delta_n   &= \pi^{n'}_{n'} -  \pi^n_n+ \pi^0_0 - \pi^n_n \notag \\
&= a_n \, \omega^n -  \, \im \, u \, \bar{\omega^n}
- A_n \, \theta \n
\end{align}
for coefficients  $\, a_n, \, A_n$. But $\, \Delta_t-\Delta_n$ is
purely imaginary, and comparing with \emph{Step 3}, $\, a_n = - 3 \, \im \,\bar{u}$.

\emph{Step 6. $\;$}   Now by considering $\, \theta$-terms in
\emph{Step 1, 2, 3, 4, 5} and the fact $\, \pi^a_i \wedge \omega^i +
\pi^a_{2n+1} \wedge \theta = 0$,  we obtain the following simple
structure equations. We omit the details of computations.
\begin{align}
\begin{pmatrix}
\pi^{n'}_{s'} - 2 \pi^n_s \n \\
\pi^{n'}_{2n+1}
\end{pmatrix}
&=
\begin{pmatrix}
0  & u \\
u  & 0
\end{pmatrix}
\begin{pmatrix}
\omega^n \\
-\, \mbox{i} \, \bar{\omega}^s
\end{pmatrix}, \\
\begin{pmatrix}
\pi^{t'}_{s'} -  \pi^t_s \\
-\pi^n_s \\
\pi^{t'}_{2n+1}
\end{pmatrix}
&=
\begin{pmatrix}
0\\
0\\
0
\end{pmatrix}, \n \\
\pi^t_{2n+1} &= (A - \,  \mbox{i} \, u \, \bar{u} ) \, \omega^t + B_t \, \theta \n \\
\pi^n_{2n+1} &= A \, \omega^n + B_n \, \theta \n
\end{align}

\emph{Step 7. $\;$}  Differentiating  $\, \pi^{t'}_{s'}- \pi^t_s =0,
\, \pi^{n}_s=0$, we get first $\, B_s =0, \, B = 0$, and
\begin{equation}\label{155}
A - \bar{A} = \im \, ( u \, \bar{u} -1 ).
\end{equation}

\emph{Step 8. $\;$}  Differentiating  $\, \pi^{n'}_n = 2 \, \omega^n
+ u \, \theta$ and $\, \pi^{n'}_{2n+1} = u \, \omega^n$,
\begin{equation}\label{156}
du = u \, ( -\pi^n_n + \pi^{2n+1}_{2n+1} ) + 2 \, ( A - A_n ) \, \omega^n
+ \im \, u \, ( \,  3 \, \bar{u}  \, \omega^n  +  u \, \bar{\omega^n} )
+ u \, ( 2 \, A - A_n) \, \theta.
\end{equation}

\emph{Step 9. $\;$}  Differentiating  $\, \pi^{n'}_{s'} = - \im \,
u \, \bar{\omega}^s$ using (\ref{156}) and collecting terms in $\,
\theta \wedge \bar{\omega}^s$,
\begin{equation}
A_n = 2 \, A - \, \bar{A}.\n
\end{equation}

\emph{Step 10. $\;$}   Differentiating $\, \pi^t_{2n+1} = (\,
\bar{A} - \im \, ) \, \omega^t$, $\, \pi^n_{2n+1} = A \, \omega^n$,
we get
\begin{equation}\label{158}
dA= \pi^0_{2n+1} + (\, A + \im \,) \, (\, \pi^{2n+1}_{2n+1} - \pi^0_0 \, )
- (\, A + \im \, )^2 \, \theta.
\end{equation}
We normalize $\, A = \im \, \alpha$ for  a real number $\, \alpha$
using group action by $\, \pi^0_{2n+1}$. Since $\, \pi^i_i + \bar{\pi}^i_i=0$,
$\, \Delta_t + \bar{\Delta}_t = \pi^0_0 + \bar{\pi}^0_0$ ana (\ref{158}) is now
reduced to
\begin{align}
d\alpha&= 2\, \im \, (\alpha+1) \,(\bar{u} \, \omega^n - u \, \bar{\omega}^n) \label{159}\\
\pi^0_{2n+1}&= -(\alpha+1)^2 \, \theta.\n
\end{align}
When $\, u \ne 0$, we may also rotate $\, u$ to be a positive
number, in which case it is determined by (\ref{155})
\begin{equation}\label{433}
2 \, \alpha + 1 = u \, \bar{u}.
\end{equation}

At this stage, note that  the only independent coefficients in the
structure equations are $\, \alpha, \, u$, and that the expression
for their derivatives does not involve any new variables. The
structure equations for the CR-flat submanifold $\, M^n \subset
\Sigma^{2n}$ thus \emph{close up} as follows.

\begin{align}
\begin{pmatrix}
\pi^{p'}_q  & \pi^{p'}_n \\
\pi^{n'}_q  & \pi^{n'}_n
\end{pmatrix}
&=
\begin{pmatrix}
\delta_{pq} \, \omega^n  &  \omega^p \\
0  & 2 \, \omega^n + \, u \, \theta
\end{pmatrix}\label{162} \\
\begin{pmatrix}
\pi^{n'}_{s'} - 2 \pi^n_s \\
\pi^{n'}_{2n+1}
\end{pmatrix}
&=
\begin{pmatrix}
0  & u \\
u  & 0
\end{pmatrix}
\begin{pmatrix}
\omega^n \\
-\im \, \bar{\omega}^s
\end{pmatrix}\n \\
\begin{pmatrix}
\pi^{t'}_{s'} -  \pi^t_s \\
-\pi^n_s \\
\pi^{t'}_{2n+1}
\end{pmatrix}
&=
\begin{pmatrix}
0\\
0\\
0
\end{pmatrix} \n
\end{align}
\begin{align}
\pi^t_{2n+1} &= (A - \im \, u \, \bar{u} ) \, \omega^t  \n \\
\pi^n_{2n+1} &= A \, \omega^n \, \n \\
\Delta_t  &=  - 2 \,  \im  \,\bar{u} \, \omega^n  - A \, \theta  \n  \\
\Delta_n  &= -3 \,  \im  \, \bar{u}  \, \omega^n -  \,  \im  \, u \, \bar{\omega^n} \,
- A_n \, \theta \n  
\end{align}
\begin{align}
du &= u \, ( -\pi^n_n + \pi^{2n+1}_{2n+1} )
+  \im  \, u \, ( \,  3 \, \bar{u}  \, \omega^n  +  u \, \bar{\omega^n} )  \n \\
&\; + 2 \, ( A - A_n ) \, \omega^n + u \, ( 2 \, A - A_n) \, \theta \n  \\
d\alpha&= 2\,  \im  \, (\alpha+1) \,(\bar{u} \, \omega^n - u \, \bar{\omega}^n)  \n \\
\pi^0_{2n+1}&= -(\alpha+1)^2 \, \theta \n
\end{align}
where $\, A= \im \, \alpha, \, A_n= 3\, \im \, \alpha$, and
$\, \alpha, \, u$ satisfy the relation (\ref{433}).
Moreover, a long but direct computation shows that these structure
equations are compatible, that is  $\, d^2 = 0$ is a formal identity of
the structure equations.

\emph{Remark.}$\;$ The structure equation closes up at order 3 [Ha].
Since both $\, \Sigma^n$ and $\, \Sigma^{2n}$ are real analytic,
this implies a $\, C^3$ nonlinear CR-flat submanifold
$\, f:\, \Sigma^n \hookrightarrow \Sigma^{2n}$ is real analytic.

Let $\, \Sigma^* = \, \Sigma^n - \, \Sigma^n_s$, and note that it is
a connected set. Note also the structure equation (\ref{162})
implies that the set of points where $\, u = 0$ or equivalently $\, \alpha = -
\frac{1}{2}$ cannot have any interior on $\, \Sigma^*$. We claim the
invariant $\, \alpha$ takes any value $\, > - \frac{1}{2}$ on $\,
\Sigma^*$.

Suppose $\, \alpha_{+} = \sup_{\Sigma^*} \, \alpha  \, > -
\frac{1}{2}$  is finite. Applying the existence part of  Cartan's
generalization of Lie's third fundamental theorem on closed
structure equations, [Br3], there exists for any $\, p_0 \in
\Sigma^*$ a neighborhood $\, U \subset \Sigma^*$ and a CR immersion
$\, g: \, U \to \Sigma^{2n}$ with  invariant $\, \alpha|_{p_0} =
\alpha_+$, hence necessarily $\, u_{p_0} = \sqrt{2 \, \alpha_{0}+1}$.
From (\ref{155}), $\, d\alpha|_{p_0} \ne 0$ and let
$\, p_{-} \in U$ be a point with $\, - \frac{1}{2} < \alpha|_{p_{-}}
< \alpha_+$. Then by uniqueness part of Cartan's theorem, there
exists a neighborhood $\, U' \subset U$ of $\, p_{-}$ on which $\,
g$ agrees with the Whitney map $\, \Gamma_n$ up to automorphisms of
$\, \Sigma^n$ and $\, \Sigma^{2n}$. Since $\, g$ and $\, \Gamma_n$
satisfy the closed set of structure equations, they are real analytic.
Thus $\, g$ is a part of $\, \Gamma_n$. But $\, du|_{p_0} \ne 0$,
and there exists a point $\, p_{+} \in U \subset \Sigma^*$ such that
$\,  \alpha|_{p_{+}} > \alpha_+$, a contradiction. By similar
argument, $\, \inf_{\Sigma^*} \, \alpha  \, = - \frac{1}{2}$, and the
claim follows for $\, \Sigma^*$ is connected.

\emph{Proof of Theorem \ref{whit}.}$\;$ Since the set of points
$\, \alpha = - \frac{1}{2}$ cannot have any interior,  let
$\, p \in M$ be a point with $\, \alpha|_p > - \frac{1}{2}$. Then
from the results above, there exists a point $\, q \in \Sigma^n$ such
that $\, \alpha|_q = \alpha|_p$. By similar argument as above after
identifying $\, p \in M$ with $\, q \in \Sigma^n$, there exists
an automorphisms $\, \tau_{2n}$ of $\, \Sigma^{2n}$
such that $\, f = \tau_{2n} \circ \Gamma_n$
on a neighborhood $\, U$ of $\, p$. The theorem follows for both
$\, f$ and $\, \Gamma_n$ are real analytic. $\; \square$

\emph{Proof of Corollary \ref{whitco}.}$\;$ By the regularity
theorem [Mi], $\, F$ is real analytic up to $\, \partial F$. Since the CR
structure on $\, S^{2n+1} = \Sigma^n$ is definite, the set of points where $\,
\partial F$ has holomorphic rank $\, n$ is a dense open subset.
Since $\, F$ is not linear, there exists a point $\, p \in \Sigma^n$
where the second fundamental form does not vanish either. By
\textbf{Theorem \ref{whit}}, $\,  \partial F$ agrees with the Whitney
map $\, \Gamma_n$ in a neighborhood of $\, p$ up to automorphisms of the unit balls.
The real analyticity then implies $\, \partial F = \Gamma_n$ on
$\, \Sigma^n$, and hence $\, F = W_n$ on $\, \mathbb{B}^{n+1}$. $\; \square$

We may  apply Cartan's generalization of Lie's third fundamental
theorem and show that the Whitney submanifold provides an example of
a deformable CR-submanifold, and in particular that it is not 1-rigid.
Take a point $\, p \in \Sigma^n$ and an analytic
one parameter family of real numbers $\, \alpha_{t} > -
\frac{1}{2}$, and set $\, u_{t} = \sqrt{2 \, \alpha_{t}+1}$. Then by
the existence part of  Cartan's theorem, there exists a neighborhood
$\, U$ of $\, p$ and a one parameter family of CR immersions $\,
f_{t}: \, U \to \Sigma^{2n}$ with the induced structure equations
(\ref{162}) such that the invariants $\, \alpha, \, u$ have the
prescribed values $\, \alpha_{t}$, $\, u_{t}$ at $\, p$. Of course
this deformation is \emph{tangential} and does not actually deform
the submanifold. It is due to an intrinsic  CR symmetry of $\,
\Sigma^n$ that cannot be extended to a symmetry of the ambient
$\, \Sigma^{2n}$ along the Whitney submanifold.

\section{Proof of Example \ref{brig}}

\noindent \textbf{Example \ref{brig}}
\emph{
Let $\, M_{n,p} \hookrightarrow \mathbb{C}P^N = P(\bigwedge^n\mathbb{C}^{n+p})$
be the Pl\"ucker embedding of the Grassmannian
$\, Gr(n, \, \mathbb{C}^{n+p})$, $\, p \leq n$.
Let $\, \hat{M}_{n,p} \hookrightarrow S^{2N+1}$ be the inverse image of
$\, M_{n,p}$ under Hopf map, which is an $\, S^1$-bundle over $\, M_{n,p}$.
When $\, p = 2, \, 3$,  $\hat{M}_{n,p}$ is a Bochner rigid
CR submanifold of $\, S^{2N+1}$. }

It is easy to check that the canonical $\, S^1$-bundle $\, \hat{M} \hookrightarrow S^{2N+1}$
over any complex submanifold $\, M \hookrightarrow \mathbb{C}P^N$ is a CR submanifold.
Moreover, the fundamental forms of $\, \hat{M}$ as a CR submanifold
are simply the pull back, in an appropriate sense, of the usual fundamental forms
of $\, M$ as a projective subvariety [HY]. It thus suffices to show that
the fundamental forms of the Pl\"ucker embeddings of the complex Grassmannian
manifolds are Bochner rigid when $\, p = 2, \, 3$.
Note the height $\, \tau$ of $\, \hat{M}_{n,p}$ is  $\, p$.

Let $\, V$ denote the tangent space of $\,  M_{n,p}$,
and we follow the notations of \tb{Section 3}.

\emph{Case $\, p=2$}.$\;$
There exists a unitary basis $\, \{\,  x^i, \, y^i \, \}$, $i = 1, \, ... \, n$, of $\, V^*$
such that the second fundamental form is given by $\, \frac{n(n-1)}{2}$ quadratic forms
\[ F^{ij} = \, x^i \, y^j - \, x^j \, y^i, \quad i < j. \]
Let $\, P^{ij} \in S^{2,0}(V^*)$, $ i<j$, be such that
\[ F^{ij} \, \bar{P}^{ij} + P^{ij} \, \bar{F}^{ij}
 = (\, x^i \, \bar{x}^i + \,  y^i \, \bar{y}^i \,) \, Q    \]
for some $\,  Q \in S^{1,1}(V^*)$. Since $\, F^{ij}$ has only $\, x, \, y$ terms,
$\, Q$ cannot have $\,  x \, \bar{x}$ terms nor $\,  y \, \bar{y}$ terms, and thus
$\, Q$ only has $\,  x \, \bar{y}$, $\,  y\, \bar{x}$ terms. Consider
$\, x^1 \, \bar{x}^1 \, Q$ term. Since it contains three $\, x$'s,
the only possible contribution comes from
$\, x^1 \, y^j \, \bar{P}_{xx}^{1j} + P^{1j}_{xx} \, \bar{x}^1 \, \bar{y}^j$,
where $\, P^{1j}_{xx}$ denotes the $\, x x$ component of $\, P^{1j}$.
But
\[ \frac{\partial}{\partial \, y^1}
\, (x^1 \, y^j \, \bar{P}_{xx}^{1j} + P^{1j}_{xx} \, \bar{x}^1 \, \bar{y}^j) = 0. \]
Hence $\, \frac{\partial}{\partial \, y^1} \,Q = 0$. By permutation symmetry
in $\, x, \, y$, and in $\, i$, and since $\, Q = \bar{Q}$, $\, Q = 0$.

\emph{Case $\, p=3$}.$\;$
There exists a unitary basis $\, \{\,  x^i, \, y^i, \, z^i \,  \}$, $i = 1, \, ... \, n$,
of $\, V^*$ such that the second fundamental form is given by the quadratic forms
\begin{align}
 F_1^{ij} &= \, y^i \, z^j - \, y^j \, z^i, \quad i < j    \notag \\
 F_2^{ij} &= \, z^i \, x^j - \, z^j \, x^i, \quad i < j       \notag \\
 F_3^{ij} &= \, x^i \, y^j - \, x^j \, y^i, \quad i < j.       \notag
\end{align}
Let $\, P_1^{ij}, \, P_2^{ij}, \, P_3^{ij} \in S^{2,0}(V^*)$, $ i<j$, be such that
\[  F_a^{ij} \, \bar{P}_a^{ij} + P_a^{ij} \, \bar{F}_a^{ij}
 = (\, x^i \, \bar{x}^i + \,  y^i \, \bar{y}^i + \,  z^i \, \bar{z}^i\,) \, Q    \]
for some $\,  Q \in S^{1,1}(V^*)$. Consider this equation$\mod x, \, y$,
and  $\, z$ in turn. Since $\, Q$ is quadratic, the result for the case $\, p=2$
implies $\, Q = 0$.

The third fundamental form is given by the cubic forms
\[ F^{ijk} = x^i \, y^j \, z^k + \, x^j \, y^k \, z^i +  x^k \, y^i \, z^j
 - \,  x^i \, y^k \, z^j - \, x^j \, y^i \, z^k  -  x^k \, y^j \, z^i, \quad i < j < k. \]
Let $\, P^{ijk} \in S^{3,0}(V^*)$, $ i<j<k$, be such that
\[ F^{ijk} \, \bar{P}^{ijk} + P^{ijk} \, \bar{F}^{ijk}
 = (\, x^i \, \bar{x}^i + \,  y^i \, \bar{y}^i + \,  z^i \, \bar{z}^i\,) \, Q    \]
for some $\,  Q \in S^{2,2}(V^*)$. Considering this equation$\mod x, \, y$,
and  $\, z$ in turn, every monomials of $\, Q$ is either a $\, x x  y  z$,
$\, y  x  y   z$, or $\, z x  y z$ term(ignoring the conjugation).
Consider $\, x  x  y  z$ component $\, Q_{xxyz}$ and the
$\, x^1 \, \bar{x}^1 \,  Q_{xxyz}$ term in the above equation.
Since there are 4 $x$-terms, the only possible contribution comes from
$\,  x^1 \, y^j \,  z^k \, \bar{P}_{xxx}^{1jk} + P^{1jk}_{xxx} \,
 \bar{x}^1 \, \bar{y}^j \,  \bar{z}^k$,
where $\, P^{1jk}_{xxx}$ denotes the $\, x  x  x$ component of $\, P^{1jk}$.
But
\[ \frac{\partial}{\partial \, y^1}
\, (x^1 \, y^j \,  z^k \, \bar{P}_{xxx}^{1jk} + P^{1jk}_{xxx} \,
 \bar{x}^1 \, \bar{y}^j \,  \bar{z}^k) = 0. \]
Hence $\, \frac{\partial}{\partial \, y^1} \,Q_{xxyz} = 0$.
By permutation symmetry in $\, x, \, y, \, z$, and in $\, i$,
and since $\, Q = \bar{Q}$, $\, Q = 0$.

\vspace{2pc}

\noindent
\begin{center}
\textbf{\large{References}}
\end{center}

\noindent [BG]  Bryant, Robert L.; Griffiths, Phillip,
Characteristic cohomology of differential systems. I. General theory,
J. Amer. Math. Soc.  8  (1995),  no. 3, 507--596

\noindent [Br1]  Bryant, Robert L., \emph{Lectures on the Geometry
of Differential Equations}, to appear

\noindent [Br2]  \underline{\quad \quad}, A duality theorem for
Willmore surfaces, J. Differential Geom.  20  (1984),  no. 1, 23--53

\noindent [Br3]  \underline{\quad \quad}, Bochner-K\"ahler metrics,
J. Amer. Math. Soc.  14  (2001),  no. 3, 623--715

\noindent [BBG]  \underline{\quad \quad}; Berger, Eric; Griffiths,
Phillip, The Gauss equations and rigidity of isometric embeddings,
Duke Math. J.  50  (1983),  no. 3, 803--892

\noindent [Ch]  Chern, S. S.,
The Geometry of $\, G$-structures, S. S. Chern Selected papers vol III, 23--76

\noindent [ChH]  Cho, Chung-Ki; Han, Chong-Kyu,
Finiteness of infinitesimal deformations and infinitesimal rigidity of hypersurfaces
in real euclidean spaces,  Rocky Mt. J. Math. 33 (2005), to appear

\noindent [ChM]  Chern, S. S.; Moser, J. K., Real hypersurfaces in
complex manifolds, S. S. Chern Selected papers vol III, 209--262

\noindent [EHZ]  Ebenfelt, P.; Huang, Xiaojun; Zaitsev, Dmitri,
Rigidity of CR-immersions into Spheres, Comm. Anal. Geom.  12
(2004), no. 3, 631--670

\noindent [Fa]  Faran, James J., Maps from the two-ball to the
three-ball, Invent. Math.  68  (1982), no. 3, 441--475

\noindent [Ga]  Gardner, Robert B., \emph{The method of equivalence
and its applications}, CBMS-NSF Regional Conference Series in
Applied Mathematics, 58.  (1989)

\noindent [Gr]  Griffiths, P., On Cartan's method of Lie groups and
moving frames as applied to uniqueness and existence questions in
differential geometry, Duke Math. J.  41  (1974), 775--814

\noindent [GH]  Griffiths, Phillip; Harris, Joseph,  Algebraic
geometry and local differential geometry, Ann. Sci. E'cole Norm.
Sup. (4)  12  (1979), no. 3, 355--452

\noindent [Ha]  Han, Chong-Kyu,  Complete differential system for
the mappings of CR manifolds of nondegenerate Levi forms, Math. Ann.
309  (1997),  no. 3, 401--409

\noindent [Hu]  Huang, Xiaojun, On a linearity problem for proper
holomorphic maps between balls in complex spaces of different
dimensions, J. Differential Geom.  51  (1999),  no. 1, 13--33

\noindent [HJ]  Huang, Xiaojun; Ji, Shanyu,  Mapping $\bold B\sp n$
into $\bold B\sp {2n-1}$, Invent. Math.  145  (2001),  no. 2,
219--250

\noindent [HY]  Hwang, Jun-Muk; Yamaguchi, Keizo, Characterization
of Hermitian symmetric spaces by fundamental forms, Duke Math. J.
120  (2003),  no. 3, 621--634

\noindent [Iw]  Iwatani, Teruo, K\"ahler submanifolds with vanishing
Bochner curvature tensor, Mem. Fac. Sci. Kyushu Univ. Ser. A 30
(1976), no. 2, 319--321

\noindent [Ja]  Jacobowitz, Howard,
An introduction to CR structures.
Mathematical Surveys and Monographs 32,  AMS (1990)

\noindent [KT]  Kaneda, Eiji; Tanaka, Noboru, Rigidity for isometric
imbeddings, J. Math. Kyoto Univ.  18  (1978), no. 1, 1--70

\noindent [La]  Landsberg, J. M., On the infinitesimal rigidity of homogeneous varieties,
Compositio Math.  118  (1999),  no. 2, 189--201

\noindent [LM]  \underline{\quad \quad}; Manivel, Laurent,
On the projective geometry of rational homogeneous varieties,
Comment. Math. Helv.  78  (2003),  no. 1, 65--100

\noindent [Mi]  Mir, Nordine, Analytic regularity of CR maps into
spheres, Math. Res. Lett.  10  (2003),  no. 4, 447--457

\vspace{3pc}

\noindent

\noindent Sung Ho Wang \\
\noindent Department of Mathematics \\
\noindent Kias \\
\noindent Seoul, Corea 130-722 \\
\texttt{shw@kias.re.kr}

\end{document}